\documentclass[12pt,oneside]{article}
\newcommand{\Path}{}
\usepackage{ \Path Paper_style}
\usepackage{ \Path Keywords}
\usepackage{ \Path Environments}
\newcommand{\figs}{}
\usepackage{tikz-cd}
\usepackage{tikzcd2}
\usepackage{tabularx} \newcolumntype{L}{>{\raggedright\arraybackslash}X}

\DeclareMathOperator{\Corel}{Correl}
\DeclareMathOperator{\AutoCorel}{Auto}
\DeclareMathOperator{\Cylndr}{Cyl}
\newcommand{\cell}{\mathfrak{j}}
\renewcommand{\InvMes}{\text{Inv-Meas}}
\renewcommand{\InvProb}{\text{Inv-Prob}}
\begin{document}
	\title{Reconstructing dynamical systems as zero-noise limits}
	\author{Suddhasattwa Das\footnotemark[1]}
	\footnotetext[1]{Department of Mathematics and Statistics, Texas Tech University, Texas, USA}
	\date{\today}
	\maketitle
	
	\begin{abstract} 
		A dynamical system may be defined by a simple transition law - such as a map or a vector field. The objective of most learning techniques is to reconstruct this dynamic transition law. This is a major shortcoming, as most dynamic properties of interest are asymptotic properties such as an attractor or invariant measure. Thus approximating the dynamical law may not be sufficient to approximate these asymptotic properties. This article presents a method of representing a discrete-time deterministic dynamical system as the zero-noise limit of a Markov process. The Markov process approximation is completely data-driven. Besides proving a low-noise approximation of the dynamics the process also approximates the invariant set, via the support of its stationary measures. Thus invariant sets of arbitrary dynamical systems, even with complicated non-smooth topology, can be approximated by this technique. Under further assumptions, we show that the technique performs a convergent statistical approximation as well as approximations of true orbits.
	\end{abstract}
	
	\section{Introduction} \label{sec:intro}

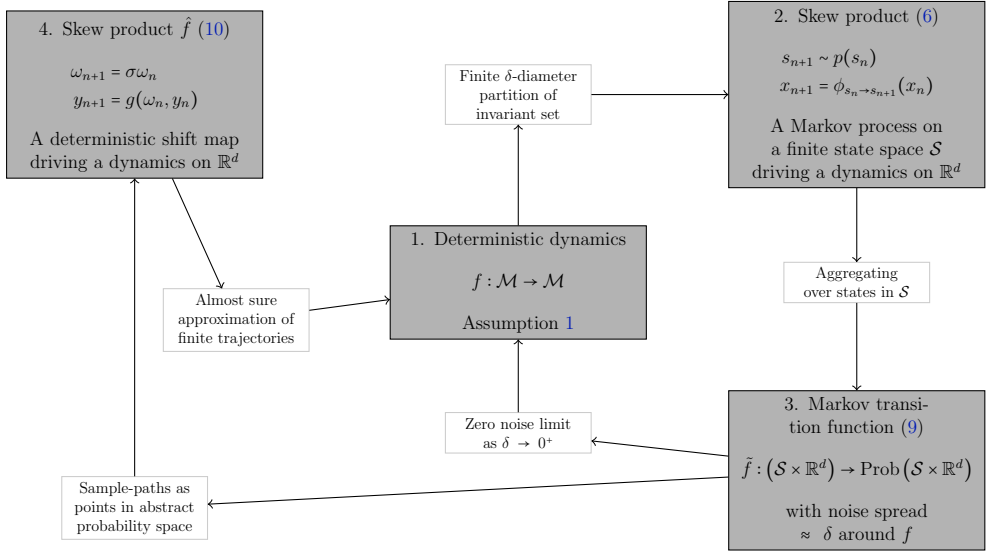
\begin{figure}\center
\begin{tikzpicture}[scale=0.5, transform shape]
	\node [style={rect11}, scale=1.2] (S1) at (0.5\columnA, -0\rowA) { 1. Deterministic dynamics \[f : \calM \to \calM \] Assumption \ref{A:f} };
	\node [style={rect11}, scale=1.2] (S2) at (2\columnA, 2.5\rowA) { 2. Skew product \eqref{eqn:Mrkv:1} \[\begin{split} s_{n+1} &\sim p(s_n) \\ x_{n+1} &= \phi_{ s_n\to s_{n+1} } ( x_n ) \end{split}\] A Markov process on a finite state space $\calS$ driving a dynamics on $\real^d$};
	\node [style={rect11}, scale=1.2] (S3) at (2\columnA, -2.5\rowA) { 3. Markov transition function \eqref{eqn:Mrkv:2} \[ \tilde{f} : \paran{ \calS \times \real^d } \to \Prob \paran{ \calS \times \real^d } \] with noise spread $\approx \delta$ around $f$};
	\node [style={rect11}, scale=1.2] (S4) at (-1.2\columnA, 2.5\rowA) { 4. Skew product $\hat{f}$ \eqref{eqn:Mrkv:3} \[\begin{split} \omega_{n+1} &= \sigma \omega_n \\ y_{n+1} &= g( \omega_n, y_n) \end{split}\] A deterministic shift map driving a dynamics on $\real^d$ };
	\node [style={rect2}] (S1_to_S2) at (0.5\columnA, 2.5\rowA) { Finite $\delta$-diameter partition of invariant set };
	\node [style={rect2}] (S2_to_S3) at (2\columnA, 0\rowA) { Aggregating over states in $\calS$ };
	\node [style={rect2}] (S3_to_S1) at (0.5\columnA, -2.0\rowA) { Zero noise limit as $\delta\to 0^+$ };
	\node [style={rect2}] (S3_to_S4) at (-1.2\columnA, -3\rowA) { Sample-paths as points in abstract probability space };
	\node [style={rect2}] (S4_to_S1) at (-0.75\columnA, -0.5\rowA) { Almost sure approximation of finite trajectories };
	\draw[-to] (S1) to (S1_to_S2);
	\draw[-to] (S1_to_S2) to (S2);
	\draw[-to] (S2) to (S2_to_S3);
	\draw[-to] (S2_to_S3) to (S3);
	\draw[-to] (S3) to (S3_to_S1);
	\draw[-to] (S3_to_S1) to (S1);
	\draw[-to] (S3) to (S3_to_S4);
	\draw[-to] (S3_to_S4) to (S4);
	\draw[-to] (S4) to (S4_to_S1);
	\draw[-to] (S4_to_S1) to (S1);
\end{tikzpicture}
\caption{Stochastic representations of a deterministic dynamical system. A deterministic dynamical system is described by a self-map $f$ as shown in the center. If $f$ is unknown, $f$ may be reconstructed from measurements of the dynamics. This article addresses the commonly used technique of finding a partition of the phase space and studying the transitions between the cells of the partition, which is induced by $f$. This results in a Markov transition on a finite state space. The paper presents a simple idea by which this Markov process can be used to create a small-noise approximation of the deterministic dynamics of $f$. The method is theoretically simple and can also be implemented by data-driven methods. The same Markov process has two other interpretations, as shown in the figure. Each interpretation lead to a different means of approximating the original deterministic dynamics. }
\label{fig:outline}
\end{figure} 

Many physical phenomenon and human driven systems \cite[e.g.]{MustaveeEtAl_covid_2021, GiannakisDas_tracers_2019, DasEtAl2023traffic, DSSY2017_QR, DasMustAgar2023_qpd}  can be described as a dynamical system : a space $\calM$ representing a collection of states of the system, and a transformation $f: \calM \to \calM$ representing the dynamics or evolution of the states. A phenomenon may be known to be describable by such a simple mathematical object, but the precise description of the transformation law $f$ may not be known. Even a description or characterization of $\calM$ may be absent. The only information available is a dataset obtained from some measurement of the phase space $\calM$. The field of data-driven study of dynamical systems aims to develop techniques of reconstructing $(\calM,f)$ from the dataset alone.

There has been a huge body of work in this field, targeting different aspects of the dynamics and relying on different numerical techniques. All methods of reconstructing the map $f$ proceeds along the following path \cite{BerryDas_learning_2022} -- one first employs an \emph{embedding mechanism} to convert the measurement into an injective map $h : \calM \to \real^d$ into a high dimensional Euclidean space. The embedding mechanism along with the embedding itself contains some information about both $\calM$ and $f$. The learning task then becomes completing this information. In the end the target is to construct a function $F : \real^d \to \real^d$ such that the \emph{semi-conjugacy} $F \circ h = h \circ f$ holds. This means that the dynamics of $F$ in $\real^d$ mirrors the dynamics of $f$ in $\calM$. This is displayed on the left below :
\begin{equation} \label{eqn:conjug}
\begin{tikzcd}
	\calM \arrow[d, "h"'] \arrow[r, "f"] & \calM  \arrow[d, "h"] \\
	\real^d \arrow[r, "F"'] & \real^d
\end{tikzcd} , \quad 
\begin{tikzcd}
	X \arrow[d, "h"'] \arrow[r, "f"] & X  \arrow[d, "h"] \\
	\real^d \arrow[r, "F"'] & \real^d
\end{tikzcd}
\end{equation}
Such a data-driven approach does not have or provide direct access to the manifold $\calM$. Instead it helps identify a subset of the \emph{data-space} $\real^d$ which may interpreted as the image $h(\calM)$ of $\calM$. If the dynamics under $f$ is restricted to a smaller invariant subset $X$ of $\calM$, then the data and hence the output of the data-driven learning will be restricted to the smaller subset $h(X)$ of $\real^d$. This is depicted in the diagram on the right of \eqref{eqn:conjug}.

Dynamical systems are known to have multiple, even uncountably infinitely many co-existing invariant sets \cite[e.g.]{DasJim17_chaos, DSSY_Mes_QuasiP_2016, DasYorkeQuasiR2016}. If the data originates from measurements made along only one such invariant set $X$, then one loses information about the diverse phenomenon occurring elsewhere in the phase-space. Thus the outcome of data-driven techniques will always fail to yield information about the complement of $X$ in $\calM$. On the other hand, the reconstruction $F$ has its own dynamics on $\real^d$. As depicted in \eqref{eqn:conjug}, under a perfect reconstruction, $F$ is only pre-determined on the image $h(X)$. Thus $F$ may have its own collection of invariant sets in $\real^d$, separate from $h(X)$. More crucially, $h(X)$ may be an unstable set, meaning that it could be inaccessible to experiments and simulations. This is made more precise via the concept of \emph{visibility}.

Given any point $z\in \calM$, its $\omega$-limit set is the collection of accumulation points of the sequence $\braces{ f^n z }_{n=0}^{\infty}$. It is thus the collection of all those points which are approached arbitrarily closely infinitely many times by the trajectory of $z$. An $\omega$-limit set can be shown to be an invariant closed set itself. The \emph{basin} of an invariant set $X$ is the set of points $z$ whose $\omega$-limit set is contained within $X$. In other words, it is the collection of points which converge towards $X$ as the dynamics enfolds. We will call the set $X$ \emph{visible} if its basin has non-zero volume with respect to (w.r.t.) the volume measure of the ambient space $\calM$. 

The concepts of visibility and stability are strongly tied to the outcome of physical experiments. For simplicity, a physical experiment may be interpreted as a simulation of the system by randomly selecting an initial point $z$. The location of this initial point is evenly distributed with respect to some volume measure on $\calM$ which is independent from the map $f$ itself. The outcome of such an experiment will be an orbit $\braces{ f^n z }_{n=0}^{\infty}$. Being visible means that with positive probability, this orbit converges to $X$. Note that the question of probability arises because we view the selection of the initial point $z$ as a random event. Similarly, being stable means that $X$ will continue to remain visible if the initial point is slightly perturbed from $X$. Stability implies visibility, but the converse is not true. The focus of our numerical methods is the issue of visibility and stability for the reconstructed dynamics.

One can have a similar discussion in a measure theoretic context. A measure $\mu$ on $\calM$ is said to be \emph{invariant} if $\mu$ is the same as its push-forward $f_*\mu$. Given a point $a\in \calM$ let $\delta_{a}$ denote the Dirac-delta measure, which is the Borel probability measure whose support is the singleton set $\{a\}$. Note that as the point $z$ travels along its trajectory one obtains a sequence of measures $\frac{1}{N} \sum_{n=0}^{N-1} \delta_{f^n z}$. These are discrete measures called \emph{empirical measures}, as they are observation based. A point $z$ is said to be in the basin of the invariant measure $\mu$ if its empirical measures converge weakly to $\mu$. One can then similarly define the properties of stability and visibility for the measure $\mu$. 

A dynamical system may have an infinite collection $\InvMes$ invariant measures. Moreover, any convex combination of invariant measures is again an invariant measure. The sub-collection $\InvProb$ of $\InvMes$ comprising f invariant probability measures is a convex set. The terminal points of this set are called \emph{ergodic} measures. These are essentially indecomposable invariant measures of $f$. They cannot be written as a non-trivial convex combination of other invariant probability. At this point we can formally state our assumptions :

\begin{Assumption} \label{A:f}
There is a manifold $\calM$, a continuous map $f: \calM \to \calM$, an invariant ergodic measure $\mu$ of $f$ whose support $X$ is compact.
\end{Assumption}

Assumption \ref{A:f} describes the base of our setup. The origin of the data is interpreted as follows :

\begin{Assumption} \label{A:data}
There is an injective map $h: \calM \to \real^d$, an initial point $x_0$ in the basin of $\mu$, leading to the $D$-dimensional timeseries $\SetDef{ h \paran{ f^n x_0 } }{ n\in \num_0 }$. 
\end{Assumption}

Assumptions \ref{A:f} and \ref{A:data} together provide a means of interpreting any general $D$-dimensional timeseries $\braces{y_n}_{y\in \num_0}$. Our goal is to reconstruct the dynamics in such a way that the invariant set is also preserved. More precisely the goal is to find an approximation of the map $F$ from \eqref{eqn:conjug} such that  the collection of visible states w.r.t. $F$ is topologically close to $h(X)$. 

Preserving the stability under a reconstruction has been an elusive goal in most learning techniques. We present a means of approximating the dynamics as an inherently stochastic system. This will lead to a loss of the determinism of $f$, but shall guarantee an approximation of $\mu$ and $X$ which is both close and stable. This notion of stability is made more precise using Arnold's paradigm \cite{Arnold_random_1991}. 

Let $\paran{ \Omega, \Sigma }$ be a measurable space, and $\tau : \Omega \to \Omega$ be a measurable map. Now suppose there is a map $g :\Omega \times \calM \to \calM$. For each $\omega \in \Omega$, one gets a different self-map $g(\omega, \cdot)$ on $\Omega$. If the choice of $\omega$ is random, then $g$ provides a parametric description of a stochastic process on $\calM$. Now consider the following dynamical system on the product space $\Omega \times \calM$:
\begin{equation} \label{eqn:skew_stoch}
\begin{split}
	\omega_{n+1} &= \tau(\omega_n) \\
	y_{n+1} &= g\paran{ \omega_n, y_n }
\end{split}
\end{equation}
This is called a \emph{skew-product} system as the first variable evolves independently and continues to drive the dynamics in $\Omega$. Such skew-product systems \eqref{eqn:skew_stoch} provide a universal description of discrete-time stochastic dynamics as a deterministic dynamical system \cite{Arnold_random_1991}. The stochasticity in the dynamics of the $y$ variable is interpreted to originate from the randomness of the initial state $\omega_0$.

Suppose that the deterministic map $f$ from Assumption \ref{A:f} corresponds to $g(\omega_0, \cdot)$ for some point $\omega_0 \in \Omega$. Now suppose that the function $g$ depends on a third parameter $t>0$ which represents a noise-bound. Thus $g$ may be denoted as $g_t$. Let $\alpha_t$ be an invariant measure for the system corresponding to $g_t$. Then $\proj_\Omega \alpha_t$ is an invariant measure for the $\Omega$-dynamics. Suppose that as $t\to 0^+$ the projections invariant measures $\proj_\Omega \alpha_t$ converge weakly to the Dirac-delta measure $\delta_{\omega_0}$. Then $f$ is interpreted to be the \emph{zero-noise limit} (z.n.l.) of the parameterized family $\SetDef{g_t}{t>0}$ of stochastic processes. Following \cite{CowiesonYoung2005}, the ergodic system $(\Omega, \mu, f)$ is said to be \emph{stochastically stable} if for any parameterized family $\SetDef{g_t}{t>0}$, any choice of invariant measures $\alpha_t$, if $\proj_\Omega \alpha_t$ converges to $\delta_{\omega_0}$, then $\proj_{\calM} \alpha_t$ must converge weakly to $\mu$.

The projection $\proj_\Omega \alpha_t$ characterizes the spread in the parameter space $\Omega$ and thus the spread in the uncertainty on the dynamics on $\Omega$. If an ergodic system is stochastically stable, and if it is represented by a stochastic process with a small spread in uncertainty around $f$, then any invariant measure of the stochastic process must be close to $\mu$. The concept of stochastic stability provides a rigorous platform on which to assess the visibility and stability of ergodic measures. Stochastic stability is hard to be established in general, and has only been demonstrated for a special class of ergodic systems called \emph{SRB}-systems \cite{SRB_young}. 

\paragraph{Goal} Reconstruct the deterministic dynamics $f$ in Assumption \ref{A:f} as a parameterized Markov stochastic process of which $f$ is a zero-noise limit. The reconstruction procedure would be entirely data-driven, and based on Assumption \ref{A:data}.

Our plan of approximation is to construct a special instance of \eqref{eqn:skew_stoch}, called a \emph{step-skew product}. The driving dynamics will be modeled as a finite-state discrete-time Markov process. This is part of a broader effort to find alternative ways to describe a dynamical system, instead of the dynamics law alone. The transformation law induces many other phenomenon \cite[e.g.]{Das2024slice, BerryDas_learning_2022} such as invariant sets, almost Markov processes and the Koopman operator. Many of these properties are asymptotic, hence an approximation of the dynamics law alone does not guarantee an approximation of these properties. We shall prove that

\begin{corollary} \label{corr:1}
Suppose there is a dynamical system as in Assumption \ref{A:f} and an embedding as in Assumption \ref{A:data}. Then the embedded dynamics in $\real^d$ is the zero-noise limit of a stochastic dynamics on $\real^d$. Moreover, the stochastic dynamics may be chosen so that the support of any  stationary measure of the process is within (Hausdorff) distance $\delta$ of the support of $\mu$.
\end{corollary}

Corollary \ref{corr:1} is proved in Section \ref{sec:step} as a consequence of our main technical result Theorem \ref{thm:1}. A Markov process is stochastic while our targeted system is deterministic. A stochastic and a deterministic process may only be compared in terms of their invariant measures, and their statistical properties. Chaotic dynamical systems have a statistical profile which is extremely complicated \cite{Nadkarni, Das2023Koop_susp, Walters2000, DGJ_compactV_2018} and difficult to estimate reliably. Theoretical results and numerical experiments seem to suggest that a small amount of noise retains the statistical properties of a dynamical system. The loss of determinism is compensated for by an easier mathematical description \cite[e.g.]{Bowen_AxiomAflow_1975, Froyland1999ulam, froyland2001extract}. This is the essential principle of the approximation scheme. 

\begin{corollary} \label{corr:2}
Let the same assumptions as Corollary \ref{corr:1} hold. Let $N>0$ be a fixed integer. Then the Markov processes can be chosen so that with high probability a simulated trajectory of length $N$ of the Markov process \eqref{eqn:Mrkv:1} coincides with an $N$-length trajectory of the true dynamics.
\end{corollary}

Corollaries \ref{corr:2} is stated more precisely and proved as a consequence of Lemma \ref{lem:ood0}, later in Section \ref{sec:symbolic}. The z.n.l. approximation exists regardless of the nature of the ergodic system $(\Omega, f, \mu)$. In the presence of stochastic stability we shall have

\begin{corollary} \label{corr:3}
If the ergodic system $(X, \mu, f)$ is stochastically stable, then the stationary measures $\mu_\delta$ of the Markov processes \eqref{eqn:Mrkv:1} converges weakly to $\mu$ as $\delta\to 0^+$.
\end{corollary}

Corollaries \ref{corr:1}, \ref{corr:2} and \ref{corr:3} summarize which of the broad goals mentioned at the beginning of this article are achievable. The mathematical foundation and technical results are laid out in the rest of the paper.

\paragraph{Outline} We next describe our technique and its theoretical basis in Section \ref{sec:step}. There are some aspects of our technique related to symbolic dynamics, which we discuss next in Section \ref{sec:symbolic}. The details about the convergence of our scheme is discussed in Section \ref{sec:cnvrgnc}. The z.n.l. approximation also has a numerical aspect, and is based on kernel methods. Some basics of kernel integral methods is discussed in Section \ref{sec:kernel}. The numerical implementation of the theory is presented in Section \ref{sec:algo}. We next look at some numerical examples to test the technique, in Section \ref{sec:example}. The theory and its effectiveness are discussed in Section \ref{sec:conclus}. 	

\section{The technique} \label{sec:step}

Throughout this section we assume Assumptions \ref{A:f} and \ref{A:data}. We shall also assume 

\begin{Assumption} \label{A:cover}
	There is a minimal, finite measurable cover $\calU = \SetDef{U_i}{ 1\leq i \leq m}$ for $h(X)$.
\end{Assumption}

A cover is called \emph{minimal} if it has no proper sub-cover. The sets $U_i$ will be called the cells of the cover. A cover is called \emph{measurable} if its cells are measurable sets with non-zero $\mu$-measure. Open covers provide a coarse graining of the phase space or invariant space. Our goal is to keep track of the transitions between the cells to obtain an outer approximation of the dynamics. This technique has been extensively studied \cite[e.g.]{Froyland1999ulam, froyland2001extract, Hunt1998unique, vanluyten2006matrix, rubido2018entropy, adler1998symbolic, cao2022abundance} in the special case when $\calU$ is a partition. Note that such a cover creates a cover $\SetDef{ V_i := h^{-1}(U_i) }{ 1\leq i \leq m}$ for $X$ in $\calM$. 

\paragraph{Transitions} Let $h_*\mu$ denote the push-forward of the measure $\mu$ onto $\real^d$. Given any cover as above, one can compute transition probabilities $\tilde{\beta}_j$ for each cell $j$ :
\begin{equation} \label{eqn:def:tildenu}
	\tilde\beta_j \in \real^m, \quad \paran{ \tilde\beta_j }_i := \mu \paran{ f(x) \in V_i \,\rvert\, x\in V_j } = \mu \paran{ f^{-1} \paran{ V_i } \cap V_j } / \mu \paran{ V_j } .
\end{equation}
We normalize this vector of conditional probabilities based at the cell $j$, into a probability distribution 
\begin{equation} \label{eqn:def:beta}
	\beta_j \in \Prob \braces{1, \ldots, M}, \paran{ \beta_j }_i := \paran{ \tilde\beta_j }_i / \sum_{l=1}^{M} \paran{ \tilde\beta_j }_l .
\end{equation}
Our construction shall contain as a sub-system a discrete state Markov process on the set $\calS := \braces{1, \ldots, m}$. The power-set of $\calS$ is its assigned sigma-algebra. Consider the $m\times m$ matrix $\mathbb{P}$ whose $j$-th column is $\beta_j$. Let $\vec{1}_m$ denote the $m$-dimensional vector with all entries $1$. By design we have $\mathbb{P}^T \vec{1}_m = \vec{1}_m$, i.e., the matrix $\mathbb{P}$ is column-stochastic. 
Any probability measure $\alpha$ on $\calS$ is a non-negative entried $m$-dimensional vector whose entries sum to $1$. The latter condition can be stated concisely as $\bracketBig{ \vec{1}_m, \alpha } = 1$. Now note that
\[ \bracketBig{ \vec{1}_m, \mathbb{P} \alpha } = \bracketBig{ \mathbb{P}^* \vec{1}_m, \alpha } = \bracketBig{ \vec{1}_m, \alpha } = 1.\]
Thus the matrix $\mathbb{P}$ converts any probability measure $\beta$ into another probability measure. Thus this matrix can be interpreted as a Markov transition function on $\calS$. For each $j\in \calS$ we shall denote by $p(j)$ the (discrete) probability measure $\beta_j$. Note that $p(j)$ also corresponds to the $j$-th column of $\mathbb{P}$.

The Markov process on $\calS$ generated by $\mathbb{P}$ has been shown to approximate statistical properties of the original ergodic system $(\Omega, \mu, f)$ under a variety of different assumptions \cite{froyland2001extract}. A condition typically relied on is quasi-compactness \cite[se]{HuntMiller1992approx} of the transfer operator $\calP$. The latter is a property borne by the spectrum of $\calP$. If $\calP$ is quasi-compact, then it has a unique fixed point. Moreover all ergodic sums converge uniquely to the unique fixed point \cite[Thm 2.3]{DingLiZhou2002fin}. Techniques such as the simplicial triangulation method \cite{DingLiZhou2002fin} are shown to provide good convergence to the fixed stationary density. Another important property concerns the iterations of $\calP$ and is said to be of \emph{Lasota Yorke type} \cite[e.g.]{DingDuLi1993high}. Simple techniques such as continuous piecewise linear approximations are shown to be convergent \cite[see]{DingZhou2001constr}. However we avoid the need for any of these conditions as our goal is a z.n.l. approximation. A z.n.l. approximation is not necessarily a stochastic approximation. We intend to do so by constructing a Markov process not on the discrete space $\calS$ alone but on the product space $\calS \times \real^d$.

\paragraph{The Markov process} Recall that by Assumption \ref{A:data} $h$ is an injective map. Thus there is at least one function $F :\real^d \to \real^d$ such that it satisfies \eqref{eqn:conjug} when restricted to $h(X)$. We now describe a low-noise approximation of $F$. The key components in the design are the functions $\phi_{j\to i}$ for various transitions $j\to i \in \calS$. Consider the subsets 
\[ X_{j\to i} := h\paran{ V_j \cap f^{-1} \paran{V_i} } = U_j \cap F^{-1}(U_i) .\]
We define $\phi_{j\to i}$ to be any continuous extension of $F$ restricted to $X_{j\to i}$. Being extensions, they satisfy
\begin{equation} \label{eqn:def:ph_ji}
	\phi_{j\to i} : \real^d \to \overline{ U_i }, \quad \phi_{j\to i} (x) = f(x), \quad \mu-a.e. \, x\in X_{j\to i} .
\end{equation}
The functions $\phi_{j\to i}$ has two important features. Firstly, its range is confined to the closure of $U_i$. Secondly it agrees with the original map on a subset of the domain. 
We now have all the components for the Markov process - the transition probabilities $\beta_j$ \eqref{eqn:def:beta} and the maps $\phi_{j\to i}$ \eqref{eqn:def:ph_ji}. They create the following Markov process on the product space $\calS \times \real^d$ :
\begin{equation} \label{eqn:Mrkv:1}
	\boxed{ \begin{aligned}
			s_{n+1} &\sim p(s_n) \\
			x_{n+1} &= \phi_{ s_n\to s_{n+1} } ( x_n )
	\end{aligned} }
\end{equation}
The functions $\phi_{ s_n\to s_{n+1} }$ are continuous functions from $\real^d$ to $U_{s_{n+1}}$, one of the cells in the cover $\calU$. This is a skew product system in which the first set of coordinates (namely $s$) evolves autonomously, and drives the second set of coordinates. This particular type of skew products is closely related to \emph{step-skew products} which have been used to demonstrate a variety of robust and non-intuitive behavior in dynamical systems \cite[e.g.]{GorodetskiEtAl1999, KleptsynNalskii2004, IlyashenkoNegut2010, DiazGelfertRams2011rich}. To specify the Markov transition function for this system, consider the map
\begin{equation} \label{eqn:def:G}
	G : \calS \times \calS \times \real^d \to \real^d , \quad (s, s', x) \mapsto \phi_{s \to s'} (x) .
\end{equation}
This function contains the collective actions of the transition functions $\phi_{j\to i}$. Note that for each fixed $(x,s) \in \calS \times \real^d$, $G(s, \cdot, x)$ is a function from $\calS$ into $\real^d$. Thus this function pushes forward the probability measure $p(s)$ on $\calS$ into the probability measure $G(s, \cdot, x)_* p(s)$ on $\real^d$. This leads to the following Markov transition function on the space $\calS \times \real^d$ :
\begin{equation} \label{eqn:trd9e}
	Q : \calS \times \real^d \to \Prob \paran{ \calS \times \real^d }, \quad (s, x) \mapsto p(s) \times \left[ G(s, \cdot, x)_* p(s) \right] .
\end{equation}
The Markov process \eqref{eqn:Mrkv:1} thus has the transition function $Q$ in \eqref{eqn:trd9e}.

The $y$-coordinates of the Markov process takes a random walk in the neighborhood $\tilde\calU := \cup_{i=1}^m U_i$ of $h(X)$. To establish this walk as a Markov process we make the simple assumption

\begin{Assumption} \label{A:partition}
	The cover $\calU$ from Assumption \ref{A:cover} is a partition, i.e., its cells are disjoint.
\end{Assumption}

Under this assumption $\calU$ is a disjoint union of the $U_j$ and each point $x\in \tilde{\calU}$ lies within a unique cell of the partition. We show later in Lemma \ref{lem:id0l3} that the Markov transition on $\calS$ has a unique stationary measure $\nu$. The Markov process on $\tilde\calU$ has the transition function : 
\begin{equation} \label{eqn:Mrkv:2}
	\hat{f} : \tilde\calU \to \Prob \paran{ \tilde\calU }, \quad \hat{f}(x) := \sum_{i : j\to i} \delta_{ \phi_{j\to i} (x) } \mathbb{P}_{j,i} , \quad \mbox{ if } x\in U_j .
\end{equation}
Equation \eqref{eqn:Mrkv:1} presents the Markov process that approximates any dynamical system using an arbitrary cover. If that cover is a partition then the Markov process on the joint space $\calS \times \real^d$ can be interpreted as a Markov process \eqref{eqn:Mrkv:2} on $\calU$, which is a neighborhood of $h(X)$ in $\real^d$. Note that $\hat{f}$ in \eqref{eqn:Mrkv:2} assigns to each point $x$ in $\calU$ a discrete measure, i.e., a weighted sum of atomic probability measures. Equations \eqref{eqn:Mrkv:1} and \eqref{eqn:Mrkv:2} are both equivalent descriptions. The former will be the target of our numerical approximation, whereas the latter will be useful in the convergence analysis. 

Recall that  the map $F$ in \eqref{eqn:conjug} is defined uniquely only on the image set $h(X)$. In a data-driven scheme, in the absence of $\calM$, our goal is to approximate $F$ as well as $h(X)$ by confining $\calU$ to some small neighborhood of $h(X)$. 
Given any cover such as $\calU$ , its \emph{mesh-size} is the largest diameter of any of its cell $\calU$. Thus the mesh size represents the spread of the cover around the target set $h(X)$. In the following statement, given any set $A$ and $\delta>0$, $\bar{B}(A,\delta)$ shall denote the closed $\delta$-neighborhood of the set $A$.

\begin{theorem} [Markov approximation of dynamics] \label{thm:1}
	Let Assumptions \ref{A:f}, \ref{A:data} and \ref{A:cover} hold. Suppose that the mesh-size of the cover is some $\delta>0$. Then
	\begin{enumerate} [(i)]
		\item Let $\mu_\delta$ be any stationary measure of the Markov process \eqref{eqn:Mrkv:1}. Then 
		\[ \support( \mu_\delta ) \subset \bar{B} \paran{ h(X), \delta } , \quad h(X) \subset \bar{B} \paran{ \support( \mu_\delta ), 2\delta } . \]
		\item Suppose Assumption \ref{A:partition} also holds. Then the projection of the ergodic system $(\Omega, f, \mu)$ via $h$ is the zero noise limit of the stochastic dynamics \eqref{eqn:Mrkv:1}, as $\delta$ approaches zero.
	\end{enumerate}
\end{theorem} 

Theorem \ref{thm:1} is proved in Section \ref{sec:cnvrgnc}. Theorem \ref{thm:1}~(ii) thus establishes any dynamical system universally as a zero noise limit of a stochastic dynamical system. The first claim implies that the support of the $\mu_\delta$ converges in Hausdorff metric to the targeted attractor $X$. This was one one of our primary goals, creating a dynamical system whose invariant region as well as iterations are a close approximation to the original. We next look at alternate way of describing the driving dynamics, as a deterministic system instead of a Markov process.	

\section{Symbolic dynamics} \label{sec:symbolic} 

Throughout this section we assume Assumptions \ref{A:f}, \ref{A:data}, \ref{A:cover} and \ref{A:partition}. There is an index function $\cell : X \to \calS$ which assigns to every $x\in X$ the number $j = \cell (x)$ which is the index of the cell of the partition in which it lies. Then for every $x\in X$ and integer $N>0$ there is a sequence of indexes	
\[ \calU(x, N) := \paran{ \cell \circ h \circ f^0(x), \ldots, \cell \circ h \circ f^N(x) } , \]
which records the sequence of cells the point visits. We call this the \emph{ $N$-itinerary} of the point. 

\paragraph{Symbolic sequences} Given any metric space $\calA$, one can associate to it a metric space called a \emph{symbolic space}. Its points are all  infinite sequences $x_0, x_1, x_2, \ldots $ of points from $\calA$. The metric structure is given by 
\[ \dist \paran{ \braces{a_n}_{n=0}^{\infty} \,,\, \braces{a'_n}_{n=0}^{\infty} } := \sum_{n=0}^{\infty} 2^{-n} \dist_{\calA} (a_n, a'_n) . \]
We denote this metric space simply by $\calA^\omega$. This space has a natural transform on it, called the \emph{shift}-map : 
\[ \sigma : \calA^\omega \to \calA^\omega, \quad \braces{a_n}_{n=0}^{\infty} \mapsto \braces{a_{n+1}}_{n=0}^{\infty} .\]
Of special interest is the case when $\calA$ is a discrete set. In that case $\calA$ can be equipped with the Dirac-delta metric. This is the metric $\delta(x,y)$ which equals $0$ if $x=y$, and is $1$ otherwise. With this choice the metric on its symbolic space becomes 
\[ \dist \paran{ \braces{a_n}_{n=0}^{\infty} \,,\, \braces{a'_n}_{n=0}^{\infty} } := \sum_{n=0}^{\infty} 2^{-n} \delta(a_n, a'_n) . \]
Symbolic spaces are of great importance in studying the computational aspects of dynamical systems. The space $\calA$ represents a discrete valued measurement function. However it is rarely the case that all the sequences of $\calA^\omega$ are observed. One is usually interested in a smaller portion of this collection. 

\paragraph{Sub-shifts} The pair $\paran{ \calA^\omega, \sigma }$ creates a continuous dynamical system of its own. Any $\sigma$-invariant subspace of $\calA^\omega$ will be called a \emph{sub-shift}. Sub-shifts are created whenever there is a function $\phi : \calM \to \calA$ on a dynamical system such as in Assumption \ref{A:f}. In that case the collection of sequences 
\[ \SetDef{ a_n := \phi \paran{ f^n x } }{ n\in\num_0 } \]
generated by all possible choices of $x\in X$ creates a sub-shift of the total symbolic space $\calX^\omega$. 

In our case, given a finite open cover $\calU$ satisfying Assumptions \ref{A:cover}, \ref{A:partition},  the role of $\phi$ is played by the index function $\cell : X \to \calS$ which assigns to every $x\in X$ and number $j = \cell (x)$ which is the index of the cell in which it lies. Any point $x\in X$ has its own trajectory $x_n = f^n(x)$. This trajectory generates a sequence of states from $\calS$ based on the index of the cell that the trajectory currently is at, namely the sequence $\SetDef{s_n = \cell \circ h \circ f^n(x) }{n\in\num}$. Thus for every $n\in\num$, $x_n \in U_{s_n}$. This is called the \emph{itinerary} of the point $x$, with respect to (w.r.t.) the cover $\calU$. The collection of itineraries is a sub-shift of the symbolic space $\calS^{\omega}$. We call this space $\Omega$. Sub-shifts provide a discrete description of continuum dynamics, and enable concepts such as information, entropy and complexity to be associated to dynamical systems \cite[e.g.]{Adler1998, Robinson1998, BDWY2020, Das2023_CatEntropy}.

On the other hand, the Markov process on $\calS$ generated by $\mathbb{P}$ creates its own sub-shift. It is the collection of all sequences $\braces{s_n}_{n=0}^{\infty}$ such that for every $n$, $s_{n} \to s_{n+1}$ is a feasible transition according to the transition probabilities stored in $\mathbb{P}$. Equivalently, for every $n$, $\mathbb{P}_{ s_n, s_{n+1} } > 0$. We call this sub-shift $\Omega'$. A good zero noise limit would ensure that $\Omega'$ is close to $\Omega$. This comparison can be done more precisely using the topology of the symbolic space. 

\paragraph{Cylinders} Given any finite sequence $s_1, \ldots, s_N$, the \emph{cylinder} created by this sequence is the collection $\Cylndr \paran{ s_1, \ldots, s_N }$ of all points in $\calS^{\omega}$ whose first $N$  coordinates coincide with this sequence. Such a cylinder is an open set, and we declare its \emph{length} to be $N$. The collection of all such cylinders form a countable basis for $\calS^{\omega}$. Moreover, for each $N$, the $N$-length cylinders form a partition of $\calS^{\omega}$. All sub-shifts of $\calS^{\omega}$ inherit its topology and thus the cylinders form a basis for sub-shifts too. With these notations in mind \eqref{eqn:Mrkv:1} can be recast as :
\[\begin{split}
	& \sigma : \Omega' \to \Omega' \\
	& g : \Omega' \times \calD \to \calD 
\end{split}\]
where $g$ is the map
\[ g : \paran{ \braces{ s_n }_{n=0}^{\infty} , y } \mapsto \phi_{s_0 \to s_1} (y) .\]
As a result we get the map $\tilde{f}$ on the product space $\Omega' \times \calD$ depicted with a dashed arrow below : 
\begin{equation} \label{eqn:Mrkv:3}
	\begin{tikzcd} 
		\Omega' \arrow[r, "\sigma"] & \Omega' \\
		\Omega' \arrow[dr, "g"'] \arrow[u, "\proj_1"] \times \calD \arrow[dashed, Shobuj, rr, "\tilde{f}"] && \Omega' \times \calD \arrow[ul, "\proj_1"'] \arrow[dl, "\proj_2"] \\
		& \calD
	\end{tikzcd} , \quad 
	\begin{tikzcd}  \paran{ \braces{ s_n }_{n=0}^{\infty} , y } \arrow[d, mapsto] \\ \paran{ \braces{ s_{n+1} }_{n=0}^{\infty} , \phi_{s_0 \to s_1} (y) } 
	\end{tikzcd}
\end{equation}
Equations \eqref{eqn:Mrkv:1}, \eqref{eqn:Mrkv:2} and \eqref{eqn:Mrkv:3} represent equivalent descriptions of the same phenomenon. Equations \eqref{eqn:Mrkv:1} and \eqref{eqn:Mrkv:3} are skew-product systems whereas \eqref{eqn:Mrkv:2} is a Markov transition function. The skew-products \eqref{eqn:Mrkv:1} and \eqref{eqn:Mrkv:3} differ by the nature of their driving systems. For the former it is a Markov process, whereas for the latter it is a deterministic dynamical system. Equation \eqref{eqn:Mrkv:3} represents the principle of viewing random processes not as a series of random variables but as a probability space of \emph{sample paths}. Thus sub-shift $\Omega'$ is precisely the collection of all sample paths possible by the Markov process on $\calS$. Along each sample path, the variable on $\calD$ undergoes a deterministic, not-autonomous process. The measure $\bar{\nu}$ provides the prior distribution of these sample paths. The interplay of these various representations are depicted in Figure \ref{fig:outline}.

\paragraph{True trajectories} Consider the cylinder
\[ \calC(x, N) := \Cylndr \paran{ \calU(x, N) } .\]
It is therefore an open subset of all those sequences in the sub-shift $\Omega'$, whose first $N+1$ coordinates coincide with the $N$-itinerary of $x$. The main utility of these sets $\calC(x, N)$ is that they contain sample paths that faithfully recreate the deterministic trajectories of $f$. Suppose we fix an $x\in X$ and a sample path $\vec{s} = \braces{s_n}_n \in \calC(x, N)$. Then note that by \eqref{eqn:def:ph_ji} we have
\[\begin{split}
\tilde{f}^N \paran{ \vec{s}, h(x) } &:= \phi_{s_{N-1} \to s_N} \circ \cdots \circ \phi_{s_{1} \to s_2} \circ \phi_{s_{0} \to s_1} (x) = \phi_{s_{N-1} \to s_N} \circ \cdots \circ \phi_{s_{1} \to s_2} \circ F | \alpha_{s_{0} \to s_1} h(x)\\
&= \phi_{s_{N-1} \to s_N} \circ \cdots \circ \phi_{s_{1} \to s_2}  \circ h \paran{ f(x) } \\
&= \tilde{f}^{N-1} \paran{ \sigma\vec{s}, h ( f(x) ) } .
\end{split}\]
Iterating this identity $N$ times we get
\begin{equation} \label{eqn:posd9}
	\tilde{f}^N \paran{ \vec{s}, h(x) } = \tilde{f}^{0} \paran{ \sigma^{N} \vec{s}, h \paran{ f^N( x ) } } = h \paran{ f^N( x ) }, \quad \forall \vec{s} \in \calC(x, N) .
\end{equation}	
Simply put, \eqref{eqn:posd9} says that the symbolic sequences corresponding to $\calC(x,N)$ follow the true trajectory of $x$ up to time $N$. 
Suppose we randomly pick an $N+1$-length sequence $\vec{s} = \paran{ s_0, \ldots, s_N }$ from $\Omega'$. The first coordinate corresponds to the cell $U_{s_0}$ of the partition $\calU$. Next suppose that we pick randomly pick a point $x\in U_{s_0}$. We wish to estimate the probability that $\calU(x, N)$ equals $\vec{s}$. To be more precise, consider the following functions : 
\begin{equation} \label{eqn:def:theta}
	\begin{split}
		\theta( x, N ) &:= \bar\nu \paran{ \Cylndr \paran{ \cell(x) } \cap \calU(x,N) } \\
		\theta( N ) &:= \int_{x\in X} \theta( x, N ) d \mu(x)
	\end{split}
\end{equation}
Here $x$ is a point drawn from $X$. Then $\Cylndr \paran{ \cell(x) }$ is an open set of sequences in the shift space $\calS^\omega$ whose first coordinate coincides with the cell that contains $x$. The quantity $\theta( x, N )$ measures the probability of a random sequence being able to recreate the $N$-itinerary of $x$, given that its first coordinate is $\cell(x)$. Thus $\theta( N )$ is the total probability that an iteration of \eqref{eqn:Mrkv:2} from an initial state would yield a true trajectory of length $N$. We can expect this probability to diminish exponentially fast as $N$ increases. 	However for a fixed $N$ we can create the partition by a special construction.

\paragraph{Towers} One of the most fundamental discoveries in the measure theoretic study of dynamical systems was the concept of \emph{tower}s and their utility in approximating the dynamics and constructing a variety of non-intuitive behavior \cite{Kornfeld2004Rokhlin, ArnouxEtAl1985cutting, KingKalikow2018alpern, Carroll1992rokhlin, EigenPrasad1997multiple}. A \emph{Rokhlin tower} \cite{Alpern2006generic, Kalikow2012inf} with base $B$, height $N$ and error $\epsilon$ consists of a measurable set $B$ such that
\begin{enumerate} [(i)]
	\item the sets $\SetDef{ B_n := f^{n} (B) }{ 0 \leq n < N }$ are disjoint;
	\item the remainder set $B_N := X \setminus \cup_{n=0}^{\infty} B_n$ has $\mu$-measure less than $\epsilon$.
\end{enumerate} 

If we create a partition using these sets, then the only transition from cell $s$ would be to cell $s+1$ for $0\leq s<N$. Such a partition will thus create a simple Markov transition, in which the transition is mostly deterministic except in the last cell $f^{N-1} (B)$ of the Rokhlin tower. To understand the utility of this simple form we consider the \emph{symbolic dynamics} induced by the original map $f$. 

\begin{lemma} \label{lem:ood0}
	Suppose Assumption \ref{A:f} holds, and $\calU$ is a partition created out of a Rokhlin tower of height $N'+N$ and remainder $\epsilon$. Then the probability $\theta( N ) $ from \eqref{eqn:def:theta} is at least $(1-\epsilon) \frac{N'-1}{N+N'}$.
\end{lemma}

Let $B$ be the base of the tower. We shall denote the cell $f^t (B)$ by $B_t$. The key observation is that
\[ \theta( x, N ) \equiv 1, \quad  \forall 0 \leq t \leq N', \; \forall x\in B_t . \]
Using this we can draw the following estimate :
\[\begin{split}
	\theta( N ) &:= \int_{x\in X} \theta( x, N ) d \mu(x) \geq \sum_{t=0}^{N'-1} \int_{x\in B_t} \theta( x, N ) d \mu(x) = \sum_{t=0}^{N'-1} \int_{x\in B_t} 1 d \mu(x) \\
	&= \sum_{t=0}^{N'-1} \mu \paran{ B_t } = \sum_{t=0}^{N'-1}  \frac{1-\epsilon}{N+N'} = (1-\epsilon) \frac{N'-1}{N+N'}
\end{split}\]
This completes the proof of Lemma \ref{lem:ood0}. \qed  

Note that the choice of $N, N'$ and $\epsilon$ were arbitrary. We also have that
\[\lim_{\begin{array}{c} \epsilon\to 0^+ \\ N'\to \infty 	\end{array}} (1-\epsilon) \frac{N'-1}{N+N'} =1 . \]
This is the precise statement of Corollary \ref{corr:2} stated before.

\paragraph{Transformations} A partition created from a Rokhlin tower construction thus ensures a good recreation of true trajectories from a simulation of the Markov process \eqref{eqn:Mrkv:3} However the partitions are extremely irregular sets, and although they have small measure $\approx 1/N$, their diameter is big. This can be amended by a change of variables. The following lemma is a standard result from Analysis (e.g. \cite{Rokhlin1949fund},  \cite[Thm 9.3.4]{Bogachev2007meas} )

\begin{lemma} \label{lem:interval}
	Any separable, non-atomic measure space possessing a complete basis, is isomorphic to the Lebesgue measure on the unit interval.
\end{lemma}

The lemma guarantees that for each cell $\calP_t$ of the partition, there is a measure theoretic isomorphism
\[ \begin{tikzcd} \calU_t \arrow[ rr, "\Psi_t", "\cong"'] && B( a_t, \delta) \end{tikzcd} , \quad \forall \, 0 \leq t \leq N , \]
for some arbitrarily chosen points $a_t$ in $\real^d$.
The transformations $\Psi_t$ are highly irregular. Fix a constant $\eta\in (0,1)$. However by Lusin's theorem \cite[e.g.]{lusin1912meas} $\Psi_t$ coincides with a continuous function $\tilde{\Psi}_t$ on a compact subset $\tilde{\calU}_t$ of $\calU_t$ of measure at least $\eta \mu(\calU_t)$. Moreover the $\tilde{\Psi}_t$ can be chosen such that its supremum norm is less than or equal to that of $\Psi_t$.

Thus we have a collection of functions $\tilde{\Psi}_t : \tilde{\calU}_t \to \real^d$. Let $h$ be any common extension of these functions from the whole of $X$. The resulting partition on the image $h(X)$ will have mesh less than $\delta$, and a Markov transition matrix that converges to that of the original Rokhlin tower partition.

There are some continuous time dynamical systems which have the structure of a \emph{suspension-flow} \cite[e.g.]{Tucker99, Suda2022Poincare, Das2023Koop_susp}. For such systems there is usually a base set $B$ which is easy to recognize, and the rest of phase space can be interpreted as return paths to the base space. For such systems, the construction of a tower is much simpler. As discussed in this section, the partition from such a tower will have a sparse transition matrix.

\section{Convergence analysis} \label{sec:cnvrgnc}

The goal of this section is to prove Theorem \ref{thm:1}. For the sake of simplicity we shall assume throughout this section that the manifold $\calM$ is a sub-manifold of $\real^d$ so that $h$ from Assumption \ref{A:data} is an inclusion map. We shall need a few key observations. We begin by familiarizing ourselves with some theoretical background.

\paragraph{Markov transitions} Let $Y$ be a metric space, equipped with its Borel sigma-algebra. Given any Markov transition function $p : Y \to \Prob(Y)$, we define the \emph{spread} of $p$ to be
\[ \text{spread}(p) := \sup_{x\in Y} \diam \support p(x) . \]
The spread is an indication of the uncertainty associated with the Markov process. It is weaker than the more widely used notion of entropy \cite[e.g.]{Das2023_CatEntropy, BDWY2020}. However it is well suited to our Analysis. 

Note that a dynamical process is a deterministic process generated by a map $T:Y\to Y$. Such a map is also a Markov transition function, since it assigns to each point $x$ the Dirac delta measure $\delta_{Tx}$. Such a transition function has zero spread. The following lemma is an easy consequence of Analysis : 

\begin{lemma} \label{lem:el9y0x}
	Let $\calV$ be a metric space and $\SetDef{\beta_t}{t>0}$ be a family of measures such that $\lim_{t \to 0^+} \text{spread}(\beta_t) = 0$. Suppose there is a point $x$ in the intersection $\cap_{t>0} \support(\beta_t)$. Then the $\beta_t$ converge weakly to the Dirac-delta measure $\delta_{x}$.
\end{lemma}

The notion of a zero noise limit is dependent on a choice of measure $\beta$ for the driving dynamics. Lemma \ref{lem:el9y0x} provides a means of establishing a zero noise limit without a knowledge of $\beta$. It is based on metric considerations alone.	

\paragraph{Ergodicity} An invariant measure $\mu$ for a dynamical system is said to be ergodic if any invariant subset $S$ of the dynamics is either null w.r.t. $\nu$, or is of full $\nu$-measure. An ergodic measure is thus an indecomposable measure. It allowed no non-trivial partition into smaller invariant sets each carrying a fraction of the total measure. The space of invariant measures is a convex set and the ergodic probability measures are its extreme points. Ergodicity is an important part of our Assumption \ref{A:f}. It is a natural assumption to make, as any dynamical system with bounded trajectories have at least one ergodic measure \cite[see]{BogachevKrylovRockner2001, Walters2000}. In fact, the convergence of data-driven techniques can seldom be justified without the aid of ergodicity \cite{DSSY_Mes_QuasiP_2016, DasGiannakis_delay_2019, DasGiannakis_RKHS_2018, DGJ_compactV_2018}. Ergodicity leads to the following algebraic consequence :

\begin{lemma} \label{lem:id0l3}
	If the measure $\mu$ is ergodic then the transition matrix $\mathbb{P}$ of the partition is transitive, i.e., it has $1$ as a simple eigenvalue.
\end{lemma}

Note that Lemma \ref{lem:id0l3} holds irrespective of the choice of $\delta$ or cover $\calU$. The only requirement is ergodic. The matrix $\mathbb{P}$ is non-negative entried and column-stochastic. Moreover $\mathbb{P}$ is irreducible, i.e. it has no invariant subspace spanned by a proper subset of the coordinate vectors. To see why, let the contrary be assumed. So there are indices $i_1 < \ldots < i_k$ in $\braces{1, \ldots, m}$ such that $\mathbb{P}$ preserved the subspace spanned by the coordinate vectors corresponding to these vectors. This would mean that the set $\tilde{U} := \cup_{j=1}^{k} U_{i_j}$ is an invariant set of $f$. Note that $X \cap \tilde{U}$ is a proper subset of $X$ by the minimality of $\calU$. This violates the ergodicity of $\mu$. Thus there is no such invariant subspace. 

Thus $\mathbb{P}$ is irreducible. From the theory of non-negative matrices \cite[e.g.]{frobenius1912matrizen, BapatRaghavan1997nonneg} it now follows that the eigenvalue of maximum magnitude of $\mathbb{P}$ is $1$, and $1$ is a simple eigenvalue. This completes the proof of Lemma \ref{lem:id0l3}. The vector $\paran{\frac{1}{m} , \ldots, \frac{1}{m} }$ is a left eigenvector of $\mathbb{P}$ corresponding to eigenvalue $1$. It is also the unique stationary probability measure of the Markov process represented by $\mathbb{P}$. 

\paragraph{Spread of the Markov process} We next establish such a metric bound for the Markov process formulation \eqref{eqn:Mrkv:2}

\begin{lemma} \label{lem:jdj9p}
	The Markov process \eqref{eqn:Mrkv:2} has spread less than $4\delta(1 + \|f\|_{Lip} )$.
\end{lemma} 

Fix any function $x\in x$ and an index $j$ such that $x \in U_j$. Consider any transition $j\to i$ from $j$. There must a point $x'\in U_j$ such that $f(x) \in U_i$. Then note that
\[\begin{split}
	\dist \paran{ f(x), \phi_{j\to i} (x) } &\leq \dist \paran{ f(x), f(x') } + \dist \paran{ f(x'), \phi_{j\to i} (x) } \leq \|f\|_{Lip} \dist(x,x') + \diam( U_j ) \\
	&\leq \|f\|_{Lip} (2\delta) + 2\delta = 2\delta(1 + \|f\|_{Lip} ) .
\end{split}\]
Thus every point in the support of the probability $\hat{f}(x)$ is within distance $2\delta(1 + \|f\|_{Lip} )$ of $f(x)$. Thus the total spread of $\hat{f}(x)$ is not more than $4\delta(1 + \|f\|_{Lip} )$. This completes the proof of Lemma \ref{lem:jdj9p}. \qed 

\paragraph{Proof of Theorem \ref{thm:1} (ii)} Lemmas \ref{lem:el9y0x} and \ref{lem:jdj9p} together imply that $(\Omega, \mu, f)$ is indeed a zero-noise limit of the Markov process \eqref{eqn:Mrkv:1}. This completes the proof of the second claim of Theorem \ref{thm:1}. We now prove the first claim.

\paragraph{Proof of Theorem \ref{thm:1} (i)} Note that since the domain of each $\phi_{j\to i}$ is $U_i$, the supports of any stationary measure $\mu_\delta$ must be within $\tilde\calU$, which is within a $\delta$-neighborhood of $X$. This proves the first half of Theorem \ref{thm:1}~(i). The proof of the second half is more technical and depends on some further observations.

\paragraph{Stationary measure} A Markov process such as \eqref{eqn:Mrkv:1} can have multiple stationary probability measures. Each of these stationary measures can be simplified structurally : 

\begin{lemma} \label{lem:c3d04}
	Let $\tilde{\beta}$ be any stationary measure of the Markov process \eqref{eqn:Mrkv:1}. Then the projection of $\tilde{\beta}$ to $\calS$ must be a stationary measure of the process generated by $\mathbb{P}$.
\end{lemma}

Thus any stationary measure $\tilde{\beta}$ for the process \eqref{eqn:Mrkv:1} must project along the $\calS$ coordinates into $\beta$, the unique stationary measure on $\calS$. Thus any integral of any continuous function $\zeta$ with respect to $\tilde{\beta}$ can be rewritten as
\begin{equation} \label{eqn:kld0p}
	\int_{\calS \times \real^d} \zeta d \tilde{\beta} = \int_{s\in \calS} \left[ \int_{\real^d} \zeta(s,x) d\tilde{\beta}(x|s) \right] d\beta(s) .
\end{equation}
Here $\tilde{\beta}(\cdot|s)$ denotes the conditional of $\tilde{\beta}$ along the fibre $\{s\} \times \real^d$. Recall that the process \eqref{eqn:Mrkv:1} thus has the transition function $Q$ in \eqref{eqn:trd9e}. The following formula describes how the Markov process transforms any measure $\tilde{\beta}$ on $\calS \times \real^d$, stationary or not.
\begin{equation} \label{eqn:Qb}
	\begin{split}
		Q_* (\tilde{\beta}) &= \int_{ \calS \times \real^d } Q(x,s) d\tilde{\beta}(x,s) = \int_{ \calS \times \real^d } p(s) \times \left[ G(s, \cdot, x)_* p(s) \right] d\tilde{\beta}(x,s) \\
		&= \int_{\calS} p(s) \left[ \int_{\real^d} \left[ G(s, \cdot, x)_* p(s) \right] d \tilde{\beta}(x|s) \right] d\beta(s)
	\end{split}
\end{equation}
If $\tilde{\beta}$ is indeed a stationary measure, then both sides of \eqref{eqn:Qb} are equal to $\tilde{\beta}$ itself.

\paragraph{Follower-predecessor} It will be useful to establish some terminology to describe the causal relations created by the Markov process on $\calS$. Given two states $s,j\in \calS$ say that a state $s$ is a $\calS$-\emph{predecessor} of a state $j$ within $\calS$ if $\mathbb{P}_{s,j}>0$ . In other words, there is a Markov transition from $s$ to $j$ according to the process generated by $\mathbb{P}$. In that case $j$ will be called an $\calS$-follower of $s$. A state $i$ will be called an $\real^d$-follower of $s$ if there is some state $s''$ such that $ \phi_{s\to s''}^{-1} \paran{ U_{s''} \cap U_i }$ has non-zero $\tilde{\beta}(\cdot | s)$ measure. Note that if $j$ is an $\calS$-predecessor of $i$, then $i$ must be an $\real^d$-follower of $j$. Also note that the notions of predecessor and follower for $\calS$ and $\real^d$ coincide if Assumption \ref{A:cover} holds.

Now substitute $ \zeta(x,s) = \lambda(s) \xi(x) $ in \eqref{eqn:kld0p} and apply \eqref{eqn:Qb} to get
\[\begin{split}
	\int \zeta d\tilde{\beta} &= \int \zeta d \left[ Q_* (\tilde{\beta}) \right] \\
	&= \int_{s\in \calS} \int_{s'\in \calS} \lambda(s') d p(s)(s') \left[ \int_{x\in \real^d} \left[ \int_{x'\in \real^d} \xi(x') d \left[ G(s, \cdot, x)_* p(s) \right] (x') \right] d \tilde{\beta}(x|s) \right] d\beta(s) .
\end{split}\]
The inner-most integral can be simplified using the change of variables formula to get
\[ \int_{x'\in \real^d} \xi(x') d \left[ G(s, \cdot, x)_* p(s) \right] (x') = \int_{s''\in \calS} \xi \circ G(s, s'', x) d p(s)(s'') .\]
Using this substitution we get
\begin{equation} \label{eqn:pd3l3}
	\int \zeta d\tilde{\beta} = \int_{s\in \calS} \int_{x\in \real^d} \left[ \int_{s'\in \calS} \lambda(s') \left[ \int_{s''\in \calS} \xi \circ G(s, s'', x) d p(s)(s'') \right] d p(s)(s') \right] d \tilde{\beta}(x|s) d\beta(s) .
\end{equation}
Now take $\lambda = 1_{j}$ and $\xi = 1_{U_i}$. Then \eqref{eqn:pd3l3} becomes :
\[\begin{split}
	\tilde{\beta} \paran{ \{j\} \times U_i  } &= \int 1_{j} 1_{U_i} d \left[ Q_* (\tilde{\beta}) \right] \\
	&= \int_{s\in \calS} \int_{x\in \real^d} \left[ \int_{s'\in \calS} 1_{j}(s') \left[ \int_{s''\in \calS} 1_{U_i} \circ G(s, s'', x) d p(s)(s'') \right] d p(s)(s') \right] d \tilde{\beta}(x|s) d\beta(s) \\
	&= \int_{s\in \calS} \int_{x\in \real^d} p(s)(j) \left[ \int_{s''\in \calS} \paran{1_{U_i} \circ G} (s, s'', x) d p(s)(s'') \right] d \tilde{\beta}(x|s) d\beta(s)
\end{split}\]
The probability weight $p(s)(j)$ is just the matrix entry $\mathbb{P}_{s,j}$. The quantity $\paran{1_{U_i} \circ G} (s, s'', x)$ is $1$ if $\phi_{s\to s''}(x) \in U_i$, otherwise it is zero. Thus $\paran{1_{U_i} \circ G} (s, s'', x)$ contributes to the integral iff a transition from $s$ intersects $U_i$. Substituting we get
\[ \tilde{\beta} \paran{ \{j\} \times U_i  } = \sum_{s\in \calS} \mathbb{P}_{s,j} \int_{x\in \real^d} \left[ \sum_{s'' \in \calS} \mathbb{P}_{s, s''} \tilde{\beta}(\cdot | s) \paran{ \phi_{s\to s''}^{-1} \paran{ U_{s''} \cap U_i } } \right] d \tilde{\beta}(x|s) d\beta(s) . \]
All the terms in the summation are nonzero. The LHS $\tilde{\beta} \paran{ \{j\} \times U_i  }$ is non-zero iff there is an $\calS$-predecessor $s$ of $j$ such that $i$ is an $\real^d$-follower of $s$. 

In particular this means that $\tilde{\beta} \paran{ \{j\} \times U_j  }$ is non-zero for every $j$. This means that the support of $\mu_\delta$ intersects each of the cells of the cover $\calU$. Thus every point $x\in X$ is within a distance of $2\delta$ from a point in $\support(\mu_\delta)$. This completes the proof of Theorem \ref{thm:1}. \qed 

This completes the statement of all of our theoretical results. Our goal from the beginning has been to obtain a data-driven reconstruction of the invariant set $h(X)$ as well as $F$ from \eqref{eqn:conjug}. The procedure on which \eqref{eqn:Mrkv:1} is based requires three constructs -- the cover $\calU$, the transition matrix $\mathbb{P}$ and the maps $\phi_{j\to i}$. The first two can be easily computed from data. We outline some methods in Algorithms \ref{algo:1} and \ref{algo:2}. We next describe the construction of  extensions, using the tool of \emph{kernel-functions}. 

\section{Kernel methods} \label{sec:kernel} 
Any function $k: Y\times Y\to\real$ on a topological space $Y$ is a kernel function. One imposes various additional constraints so that deeper geometric properties of the underlying space $Y$ may be identified, or to generate suitable function space. We shall assume that

\begin{Assumption} \label{A:kernel}
	There is a strictly positive definite, symmetric kernel $k: \calM \times \calM \to \real$, which is $C^r$ for some $r\geq 1$, and strictly positive valued.
\end{Assumption} 

A kernel is said to be \emph{symmetric} if it is symmetric w.r.t. its two variables. A kernel $k$ is said to be \emph{positive semi-definite} if for every distinct points $x_1, \ldots, x_n \in Y $ and scalars $c_1, \ldots, c_n \in \cmplx$, $ \sum\limits_{i,j=1}^{n} c_i^* k(x_i,x_j) c_j \geq 0$. It will be called \emph{strictly positive definite} when equality holds in this equation if and only if $c_i = 0 \quad \forall i = 1,\dots,n$. The main utility of kernels in Learning theory is the fact that for each $y\in Y$, $k(y,\cdot)$ is a function on $Y$. It is called the \emph{section} of $k$ at $y$. Thus the first variable in a kernel can be viewed as a parameter to generate an family of functions. The span of the kernel sections create useful interpolation spaces. If the kernel is symmetric and s.p.d. the span is called a reproducing kernel Hilbert space or RKHS. RKHSs are collections of smooth functions which are universal approximators and also endowed with a Hilbert space structure. The span of a general kernel function can be described using the language of \emph{integral operators}.

\paragraph{Kernel integral operators} Let $\nu$ be a probability measure on $Y$. Associated with the kernel $k$ and measure $\nu$ is a kernel integral operator $K^\nu$ which transforms every function $\phi\in L^2(\mu)$ as 
\[ K^\nu \phi (x) := \int_{\calM} k(x, y) \phi(y) d\mu(y) , \]
If $\phi$ lies in $L^2(\nu)$ then so does $K^\nu \phi$. If the kernel function is continuous so is $K^\nu \phi$. If $Y$ is a manifold and if the kernel $k$ is $C^r$, then the functions in the image of the integral operator $K$ will also consist of $C^r$ functions. Kernel integral operators and associated Hilbert spaces have found widespread use in the mathematical theory of learning \cite{CuckerSmale2001}, and more recently in ergodic theory \cite{DasGiannakis_delay_2019, GiannakisDas_tracers_2019, DGJ_compactV_2018, DasGiannakis_RKHS_2018}, conditional expectation \cite{Das2023conditional, Das2024drift}, and Algebra  \cite{DasDimitris_CascadeRKHS_2019}. 

Suppose $k$ is square integrable with respect to $\nu$, and is strictly positive valued. Then $k$ can be modified by a process called \emph{Markov normalization} :
\[ k^{\nu}_{Markov} (x, y) := k(x,y) / \rho^{\nu}(x) , \quad \rho^{\nu}(x) := \int_{Y} k(x,z) d\nu(z) . \]
The function $\rho^{\nu}$ is known as the \emph{degree-vector}. Markov kernels have the property that 
\[ \int_{Y} k^{\nu}_{Markov} (x, y) d\nu(y) = 1, \quad \forall x\in \calM. \]
This allows such kernels to be interpreted as mimicking a Markov transition function. Hence the name. 

\paragraph{Range of Markov kernels} Since we assume the kernel $k$ to be strictly positive valued and $C^r$, the degree vector $\rho^\nu$ will be strictly positive valued and $C^r$ too, for every choice of $\nu$. Consider the \emph{multiplication operator}
\[ M_{\rho}^{\nu} : L^2(\nu) \to L^2(\nu) , \quad \phi \mapsto \phi \rho^\nu . \]
The operator $M_{\rho}^{\nu}$ is therefore merely pointwise multiplication by the degree vector. This operator is bounded with a bounded inverse. In fact $M_{\rho}^{\nu}$ is a symmetric, strictly positive definite bounded operator with a bounded inverse. Now consider the operator 
\[ \tilde{K}^\nu := \SqBrack{ M_{\rho}^{\nu} }^{1/2} K_{Markov}^\nu \SqBrack{ M_{\rho}^{\nu} }^{-1/2} : L^2(\nu) \to L^2(\nu) \]
Thus $\tilde{K}^\nu$ is a similarity transform of the Markov operator. Now note that
\[\begin{split}
	\paran{ \tilde{K}^\nu \phi }(x) &= \rho^{\nu}(x)^{1/2} \int k^{\nu}_{Markov} (x, y) \rho^{\nu}(y)^{-1/2} \phi(y) d\nu(y) \\
	&=  \rho^{\nu}(x)^{1/2} \int \frac{k(x,y)}{  \rho^{\nu}(x)  \rho^{\nu}(y)^{1/2} } \phi(y) d\nu(y) =  \int \frac{k(x,y)}{ \rho^{\nu}(x)^{1/2} \rho^{\nu}(y)^{1/2} } \phi(y) d\nu(y) \\
	&= \int \tilde{k}(x,y) \phi(y) d\nu(y) .
\end{split}\]
Here $\tilde{k}(x,y)$ is the symmetric kernel $\frac{k(x,y)}{ \rho^{\nu}(x)^{1/2} \rho^{\nu}(y)^{1/2} }$. In fact the kernel $\tilde{k}$ satisfies Assumption \ref{A:kernel}, see \cite[Sec 4]{DGJ_compactV_2018} for a proof. Thus the Markov normalized kernel is similar to the kernel integral transform corresponding to a kernel satisfying Assumption \ref{A:kernel}. This means among other things, that the range of $K^{\nu}_{Markov}$ is dense in the space of continous functions on the support of $\nu$.

The algorithms we develop in the next section utilize Markov kernels to construct the transition functions $\phi_{j\to i}$ as in \eqref{eqn:def:ph_ji}. The measure $\nu$ will be a discrete measure supported on a small portion of the dataset. A crucial feature of Markov kernels that we utilize is the control it provides over its range.

\paragraph{Convexity} Let $\phi$ be a continuous function on a compact set $Y$. Then the \emph{convex-hull} of $\phi$ is the collection
\[ \mbox{convex-hull} (\phi) := \SetDef{ \int \phi d \gamma }{ \gamma \in \Prob( Y ) } .\]
Let $p$ be a Markov kernel  and $P^{\nu}$ denote the integral operator corresponding to $p$ and probability measure $\nu$. Then note that
\[ \paran{P^\nu \phi}(x) = \int p(x,x')\phi(x') d\nu(x') = \int \phi p(x, \cdot) d\nu, \quad \forall \phi \in L^2(\nu), \, \forall x\in Y.  \]
The value $\paran{P^\nu \phi}(x)$ of $P^\nu \phi$ at $x$ is thus a weighted average of $\phi$. The weight is provided by the Markov kernel section $p(x,\cdot)$. This means that 
\[  \ran \paran{ P^\nu \phi } \subseteq \mbox{convex-hull} (\phi). \]
This relation is important from a numerical point of view. Typically $\nu$ represents an empirical measure borne by data. Thus reconstructing a function as a function of the form $P^\nu \phi$ ensures that the range of the function lies within the convex hull of the output values used in its training phase.

We next present the algorithmic implementation of the technique to clearly delineate the roles played by the data inputs and various algorithmic parameters. 

\section{Numerical implementation} \label{sec:algo}

\begin{table}
	\caption{Summary of parameters in the algorithms }
	\begin{tabularx}{\linewidth}{|L|L|L|L|L|}
		\hline
		Parameter & $\delta$ & $\theta$ & $\gamma$ & N \\ \hline
		Role & Mesh size of covering & Gaussian kernel & Multiplication constant to determine the Ridge regression coefficient & Number of data samples \\ \hline
		Values & $0.1$ & Selected according to \eqref{eqn:bw_select} at threshold 1\% & $0.001$ & Varied. See Table \ref{tab:param2} \\ \hline
	\end{tabularx}
	\label{tab:param1}
\end{table}	

\begin{figure}[!ht]\center
	\includegraphics[width=0.95\linewidth, height=0.5\textheight, keepaspectratio]{\figs 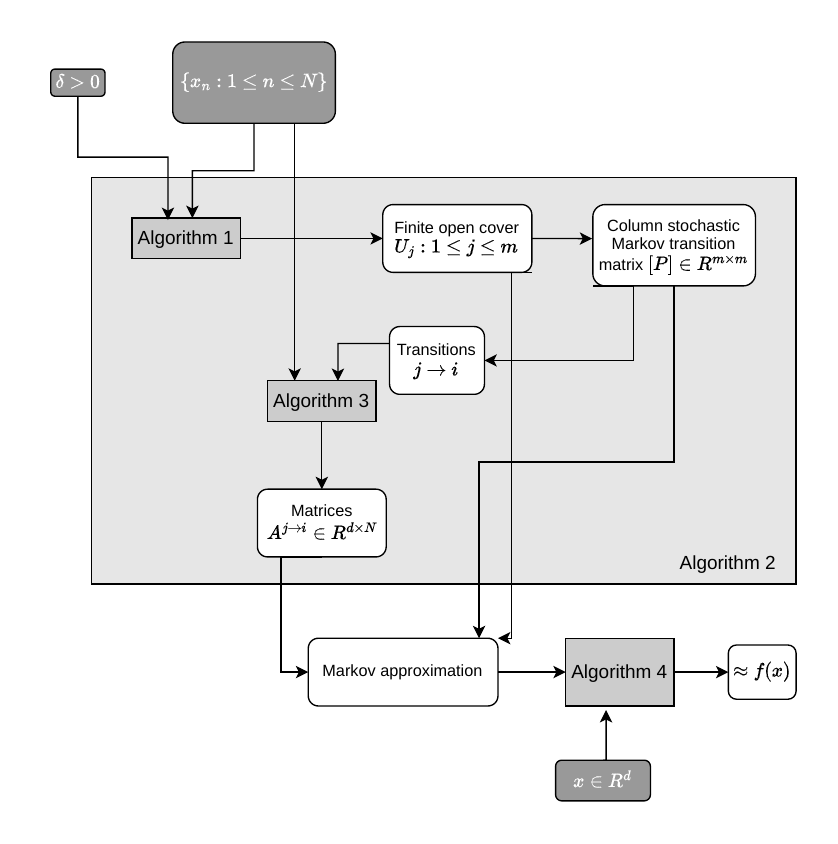}
	\caption{Outline of the components of the numerical implementation of the algorithms. The algorithms are described in Section \ref{sec:algo}.}
	\label{fig:outline1}
\end{figure}

We now describe the numerical implementation of the technique presented in Section \ref{sec:step}. The technique begins with the creation of a finite open cover. The following is an easy realization of this task.

\begin{algorithm} \label{algo:1}
	Algorithm to construct a fine-grained cover for a small neighborhood of the dataset.
	\begin{itemize}
		\item \textbf{Input.}
		\begin{enumerate}
			\item Timeseries $\calX := \SetDef{ x_n }{ 1 \leq n \leq N+1 }$ of points in $\real^d$.
			\item Grain size $\delta>0$.
		\end{enumerate}
		\item \textbf{Parameters.} 
		\item \textbf{Output.} A finite subset $J$ of $1, \ldots, N$ such that 
		\begin{equation} \label{eqn:algo:1:1}
			\cup_{j\in J} \bar{B} \paran{ \delta, x_j } \supseteq \bar{B} \paran{ \delta, \calX } = \cup_{n=1}^{N} \bar{B} \paran{ \delta, x_n } .
		\end{equation}
		\item \textbf{Steps.} Let $\chi : \real \to \real$ be the function whose value is 1 on $[0, \delta]$ and 0 outside.
		\begin{enumerate} 
			\item Initialize an empty list $J$.
			\item Construct a sparse Gaussian kernel matrix $\Matrix{K} \in \real^{N\times N}$ whose entries are $\Matrix{K}_{i, j} = \chi \paran{ \norm{ x_i - x_j } }$.
			\item Find the index $j$ corresponding to the highest row-sum.
			\item Add $j$ to $J$.
			\item Remove all vertices from the matrix within distance $\delta$ from $j$.
			\item Repeat Step 4 onwards, as long as the matrix remains non-empty.
		\end{enumerate}
	\end{itemize}
\end{algorithm}

Algorithm \ref{algo:1} is a greedy algorithm and proceeds by looking for the vertex that is within $\delta$-distance to the maximum number of remaining vertices. Algorithm \ref{algo:1} is not built to yield the optimum cover. We next extract some combinatorial information from these covers

\begin{algorithm} \label{algo:2}
	\begin{itemize}
		\item \textbf{Input.}
		\begin{enumerate}
			\item Timeseries $\calX = \SetDef{ x_n }{ 1\leq n \leq N+1 }$ of points in $\real^d$.
			\item A finite open cover $\SetDef{U_j}{ j\in \braces{1, \ldots, m} }$ for the timeseries.
		\end{enumerate}
		\item \textbf{Output.} Some combinatorial information of the transitions between the cover elements.
		\item \textbf{Steps.} For each $j\in \braces{1, \ldots, m}$, compute the following for the $j$-th cell of the cover :
		\begin{enumerate}
			\item A list $\alpha_{j \to i}$ of all cells $i$ for which there is an $n$ in $[1,N]$ such that $x_n\in U_j$ and $x_{n+1} \in U_i$.
			\item For each $i$ in $\alpha_{j \to i}$ compute the following : 
			\begin{enumerate}
				\item The collection $X_{j\to i}$ of indices $n$ in $[1,N]$ such that $x_n\in U_j$ and $x_{n+1} \in U_i$.
				\item Compute a vector $\beta_{j}$ of length $\abs{ \alpha_{j \to i} }$ such that $\paran{ \beta_{j} }_i = \abs{ X_{j\to i} }$. Then normalize $\beta_{j}$ to have total sum equal to 1.
				\item Feed the data-sets indexed by $X_{j\to i}$ and $X_{j\to i}+1$ into Algorithm \ref{algo:3} to obtain a coefficient matrix $a_{j\to i}$. 
			\end{enumerate}
		\end{enumerate}
	\end{itemize}
\end{algorithm} 

The output of Algorithm \ref{algo:2} is a list of possible transitions $j\to i$ between cells, and for each such transition, there is collection of various objects. One of the steps invokes Algorithm \ref{algo:3}, which is the following standard procedure for computing a Markov kernel from data.

\begin{algorithm} \label{algo:3}
	\begin{itemize}
		\item \textbf{Input.}
		\begin{enumerate}
			\item A collection $\SetDef{x_n}{ 1\leq n \leq M }$ of points in $\real^d$.
			\item A collection $\SetDef{y_n}{ 1\leq n \leq M }$ of points in $\real^k$.
		\end{enumerate}
		\item \textbf{Parameters.} Kernel function $k$. Ridge regression coefficient $\gamma \geq 0$
		\item \textbf{Output.} A $k\times M$ matrix $\Matrix{A}$ with columns $\paran{ a_1 , \ldots, a_M }$ such that if $\hat f = \sum_{n=1}^{M} a_n k_{Markov}(x, x_n)$ then for each $1\leq n \leq M$, $\hat f(x_n) = y_n$. 
		\item \textbf{Steps.} 
		\begin{enumerate}
			\item Compute the $M\times M$ matrix $\Matrix{K}_{i,j} := k( x_i, x_j)$ for each $1\leq i,j \leq M$.
			\item Compute the degree vector $\rho = \frac{1}{N} \Matrix{K} \vec{1}_M$.
			\item Compute the Matrix $\Matrix{P} := R^{-1} \Matrix{K}$, where $R := \diag(\rho)$.
			\item Compute $\Matrix{A}$ as the $\gamma \norm{ \Matrix{P} }$-regularized least squares solution to $\frac{1}{N} \Matrix{P} \Matrix{A} = \vec{y}$.
		\end{enumerate}
	\end{itemize}
\end{algorithm}

The object constructed in Algorithm \ref{algo:2} is the final implementation of the process \eqref{eqn:Mrkv:1}. We now show how it can be used to simulate an iteration of the Markov process.

\begin{algorithm} \label{algo:4}
	Iteration of the Markov process on an initial point.
	\begin{itemize}
		\item \textbf{Input.}
		\begin{enumerate}
			\item Training timeseries $\calX := \SetDef{ x_n }{ 1 \leq n \leq N+1 }$ of points in $\real^d$.
			\item The output tuples $\paran{ \alpha_{j \to i}, \beta_{j}, A_{j\to i} }$ of Algorithm \ref{algo:2}.
			\item Initial point $(j,x)$ with $j\in \calS$ and $x\in \real^d$.
		\end{enumerate}
		\item \textbf{Parameters.} Kernel function $k$.
		\item \textbf{Output.} A pair $(i,y)$ where $i\in \calS$ and $y\in U_i \subset \real^d$, such that the distribution of the pair depends only on $x$.
		\item \textbf{Steps.} 
		\begin{enumerate}
			\item Find a cell $U_j$ which is closest to $x$.
			\item Randomly select a transition to a cell $i$, based on the probability distribution $\beta_{j}$.
			\item Compute the pairwise distances $v_n := k(x, x_n)$ for each $n\in X_{j\to i}$.
			\item Normalize the vector $\vec{v}$ so that it sums to $1$.
			\item Compute the vector $y = A_{j\to i} \vec{v}$.
		\end{enumerate}
	\end{itemize}
\end{algorithm}

The entire numerical procedure is outlined in Figure \ref{fig:outline1}. Table \ref{tab:param1} summarizes all the parameters involved. We next look at some numerical examples to test the algorithms.

\section{Examples} \label{sec:example}

\begin{table}
	\caption{Summary of parameters in experiments. $N$ represents the number of data samples. The distance between  $L^1$-Hausdorff distance \eqref{eqn:def:Hauss_metric} }
	\begin{tabularx}{\linewidth}{|L|L|L|L|L|L|L|}
		\hline
		System & Properties & Figures & ODE Parameters & $N$ & $L^1$-Distance between support of invariant measures & Bias  \\ \hline
		Lorenz 63 \eqref{eqn:vec:l63} & Chaotic system in $\real^3$ & \ref{fig:l63} & $\sigma = 10$, $\beta = 8/3$, $\rho=28$ & $10^4$ & $0.106$ & $0.600$ \\ \hline
		Henon map \eqref{eqn:def:Henon} & Chaotic system from a discrete-time map on $\real^2$ & \ref{fig:Henon} & $a=1.4$, $b=0.3$ & $10^4$ & $0.015$ & $0.027$  \\ \hline
		Lorenz 96 \eqref{eqn:L96} & Chaotic system of $m$ 1-dimensional oscillators, cyclically coupled & \ref{fig:L96} & $f=8.0$, $m=10$ & $2\times 10^4$ & $0.338$ & $0.173$ \\ \hline
	\end{tabularx}
	\label{tab:param2}
\end{table}

\begin{figure}[!ht]\center	 
	\includegraphics[width=0.35\linewidth, height=0.5\textheight, keepaspectratio]{\figs 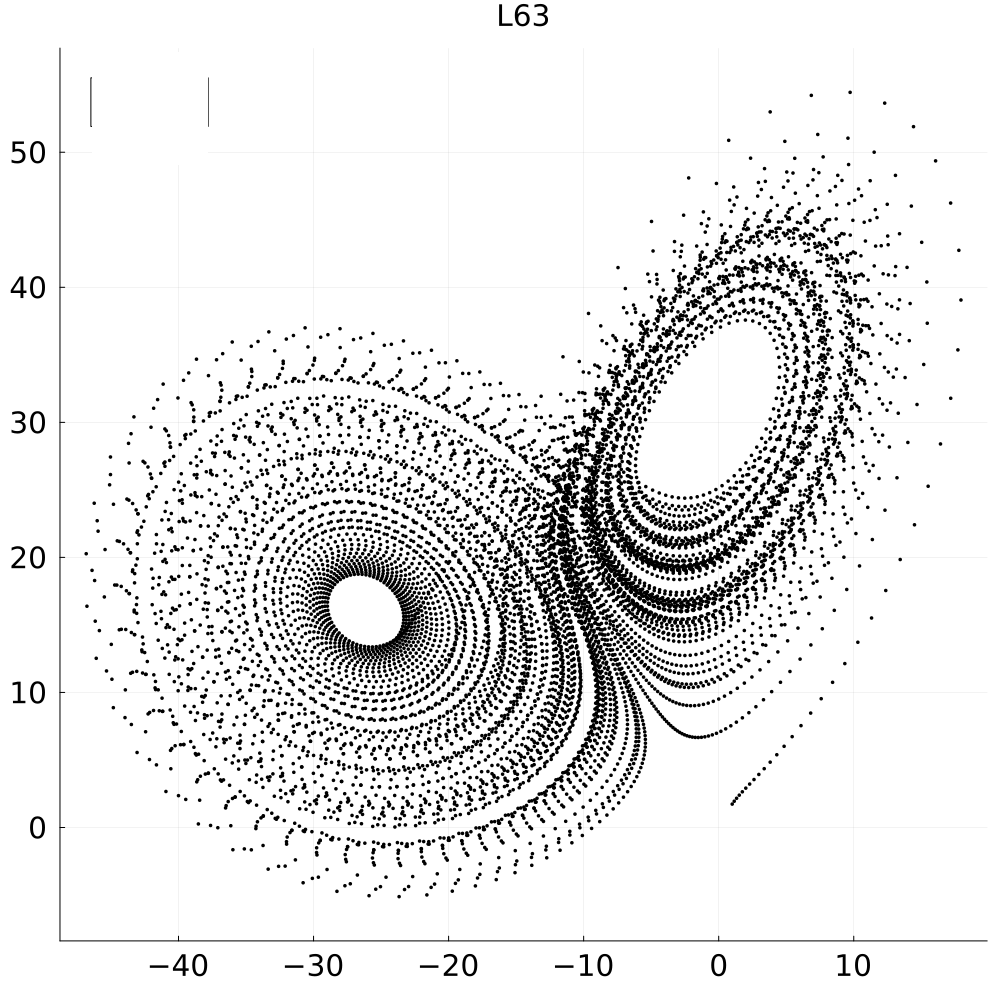}
	\includegraphics[width=0.35\linewidth, height=0.5\textheight, keepaspectratio]{\figs 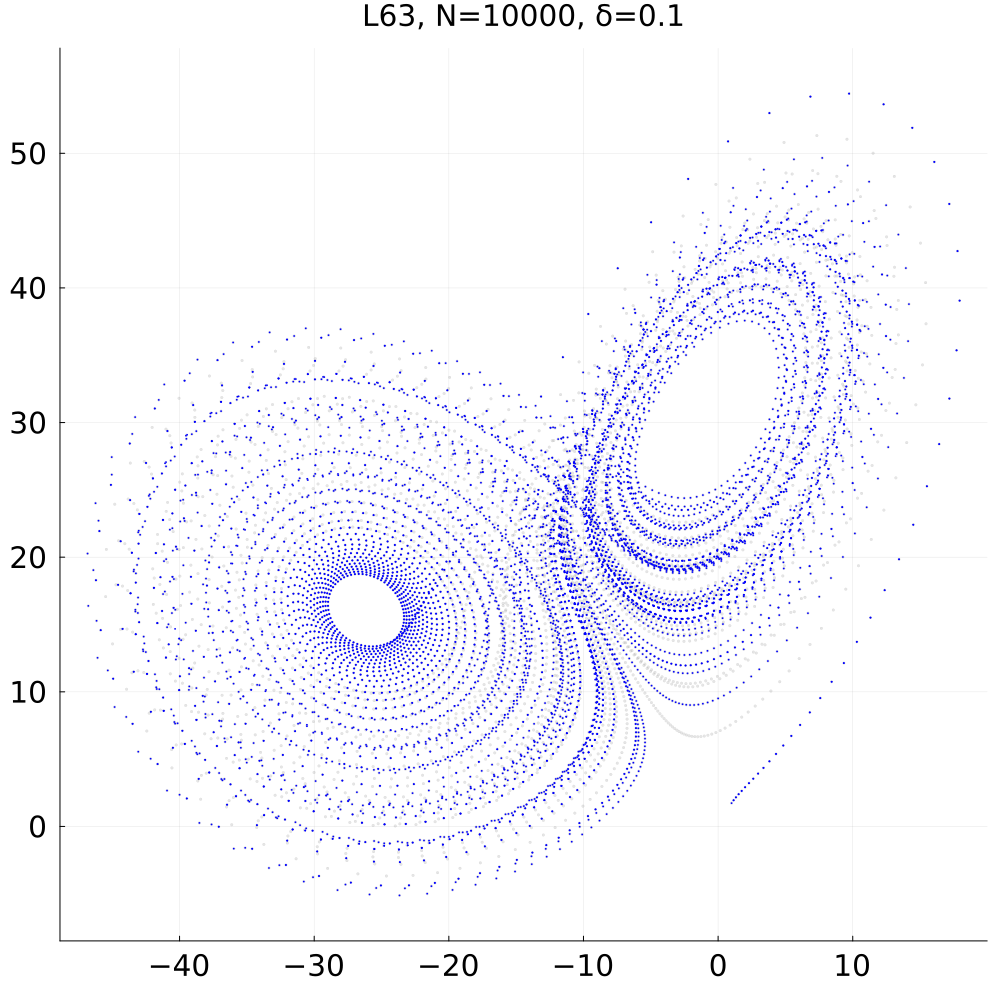}
	\includegraphics[width=0.28\linewidth, height=0.5\textheight, keepaspectratio]{\figs 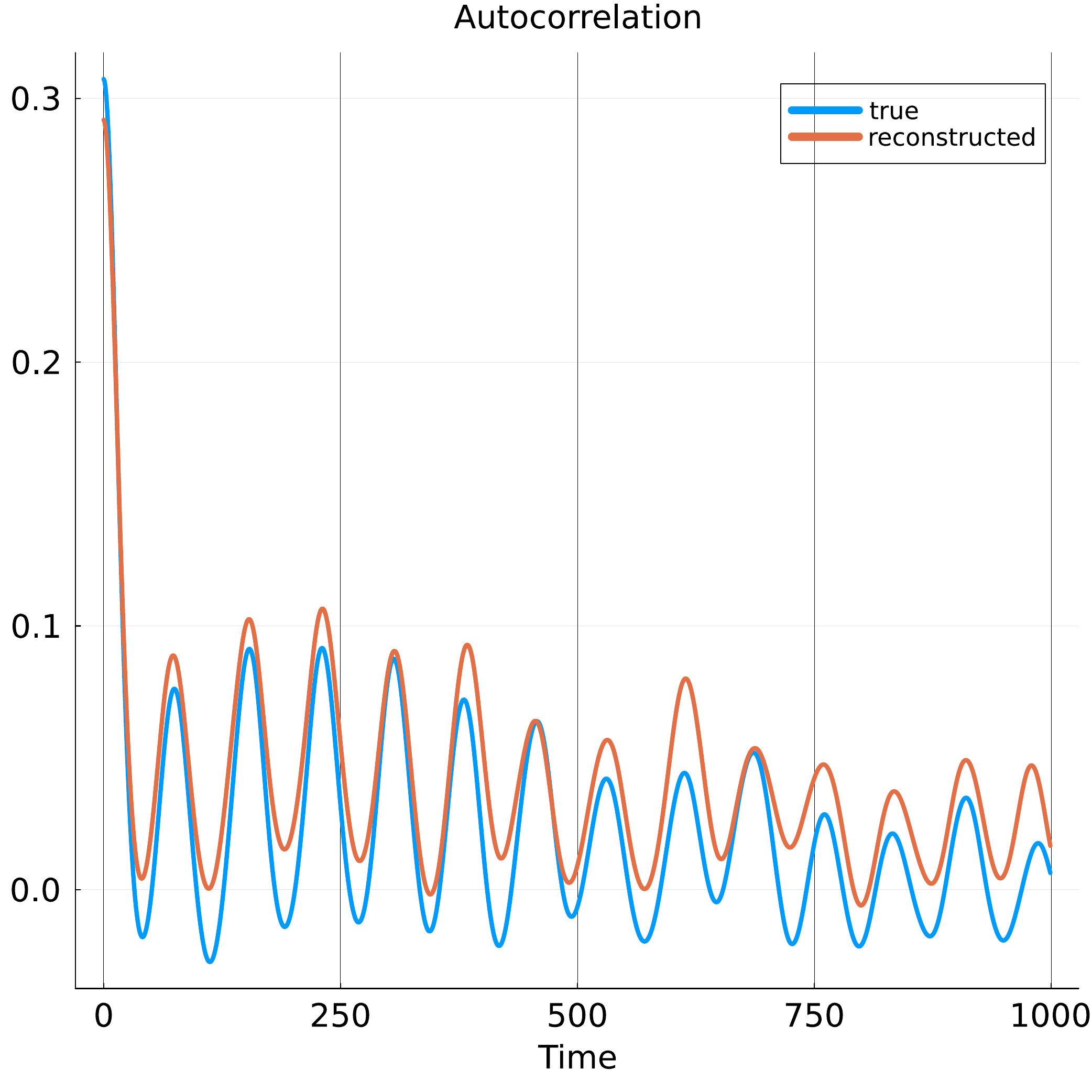}
	\caption{Reconstructing the Lorenz63 system \eqref{eqn:vec:l63}.}
	\label{fig:l63}
\end{figure}

\begin{figure}[!ht]\center	 
	\includegraphics[width=0.35\linewidth, height=0.5\textheight, keepaspectratio]{\figs 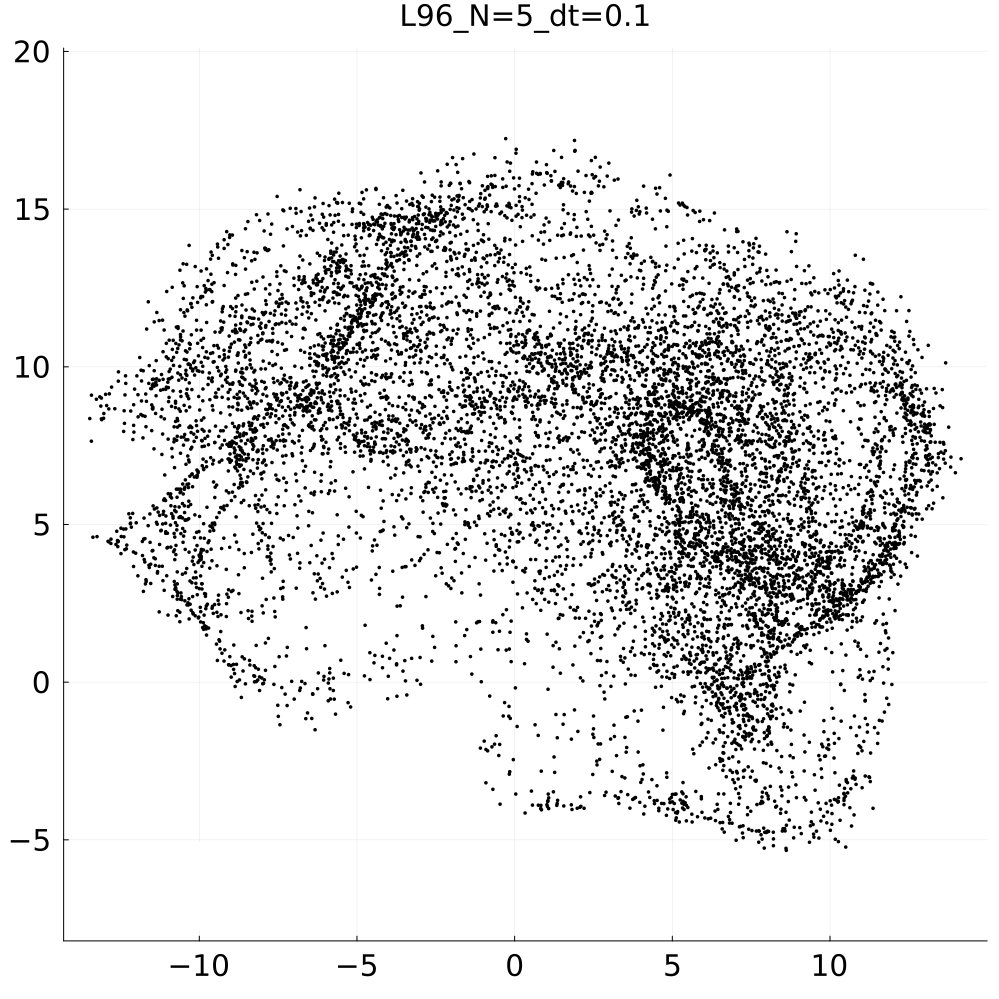}
	\includegraphics[width=0.35\linewidth, height=0.5\textheight, keepaspectratio]{\figs 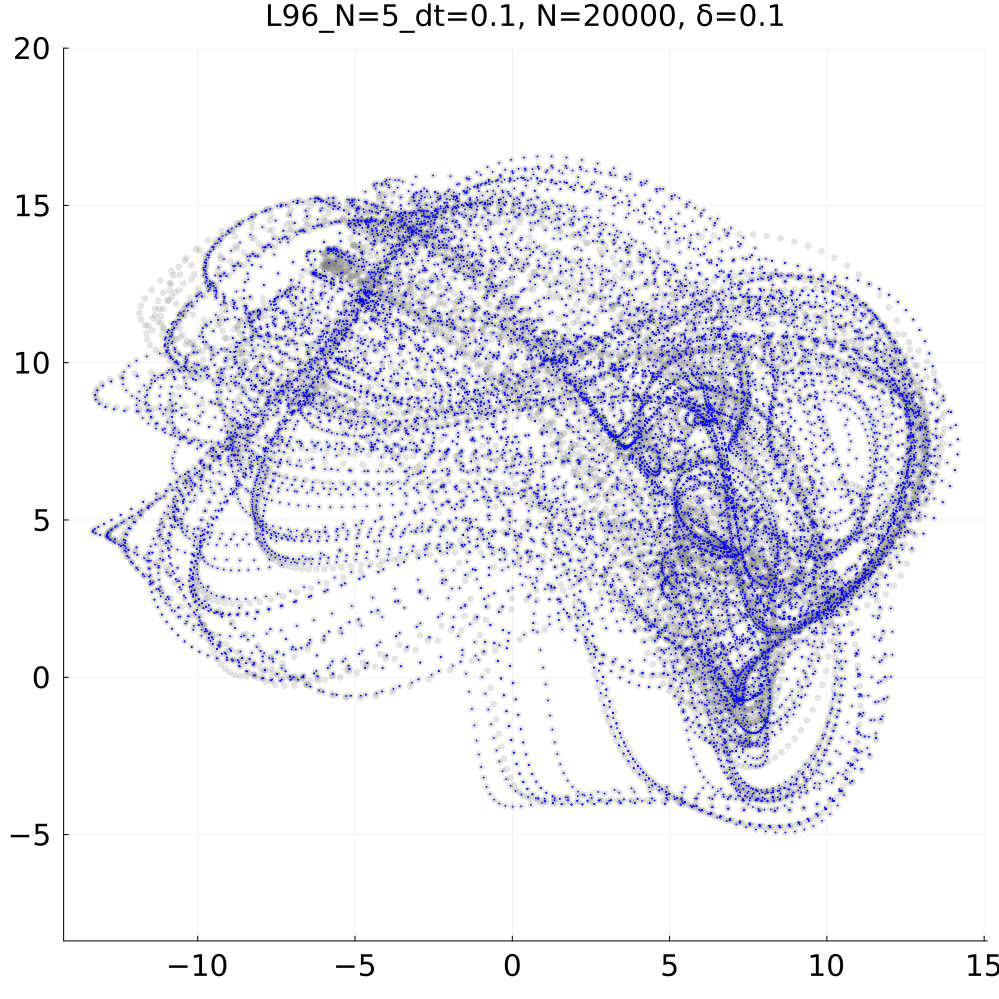}
	\includegraphics[width=0.28\linewidth, height=0.5\textheight, keepaspectratio]{\figs 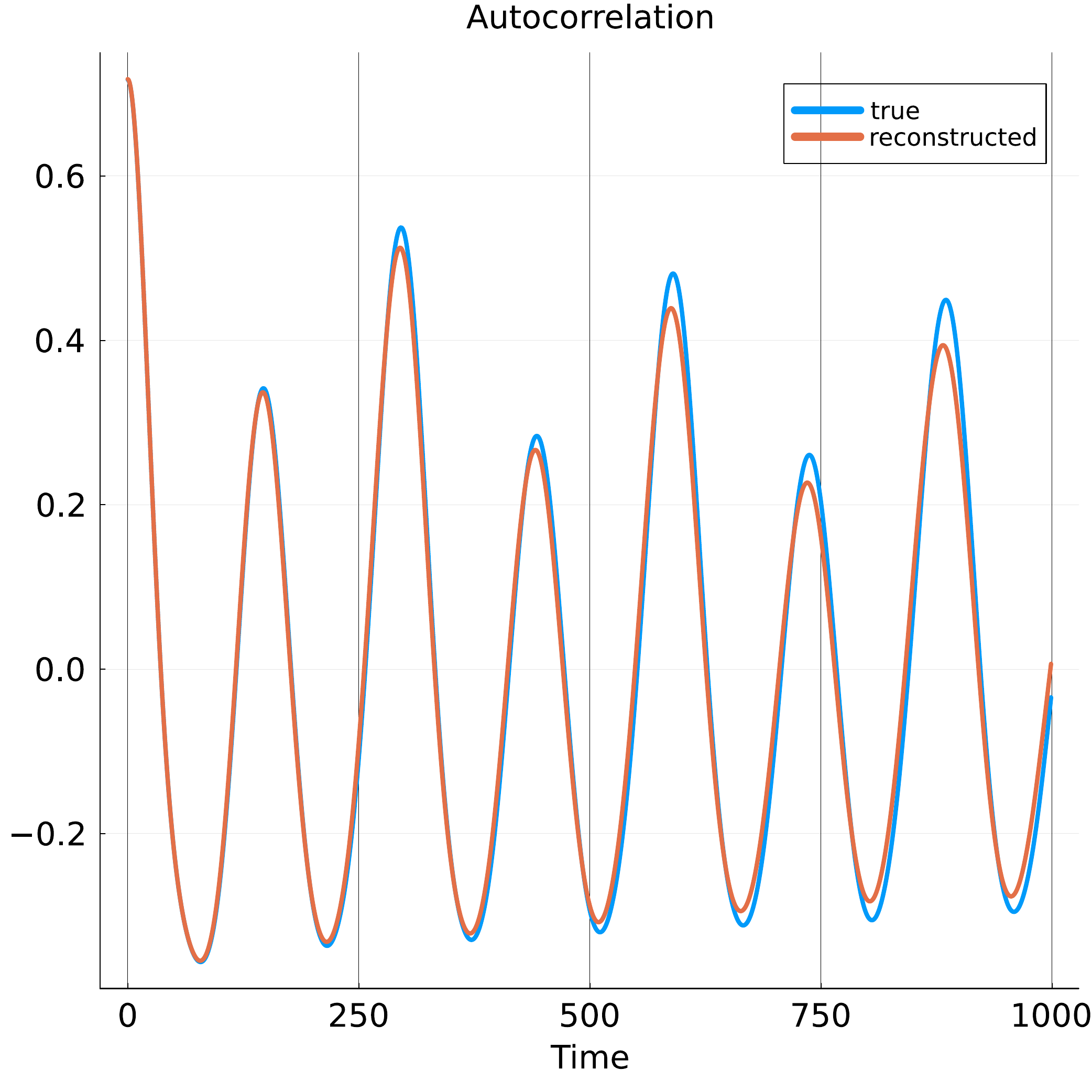}
	\caption{Reconstructing the Lorenz-96 system \eqref{eqn:L96}.}
	\label{fig:L96}
\end{figure}

\begin{figure}[!ht]\center	 
	\includegraphics[width=0.35\linewidth, height=0.5\textheight, keepaspectratio]{\figs 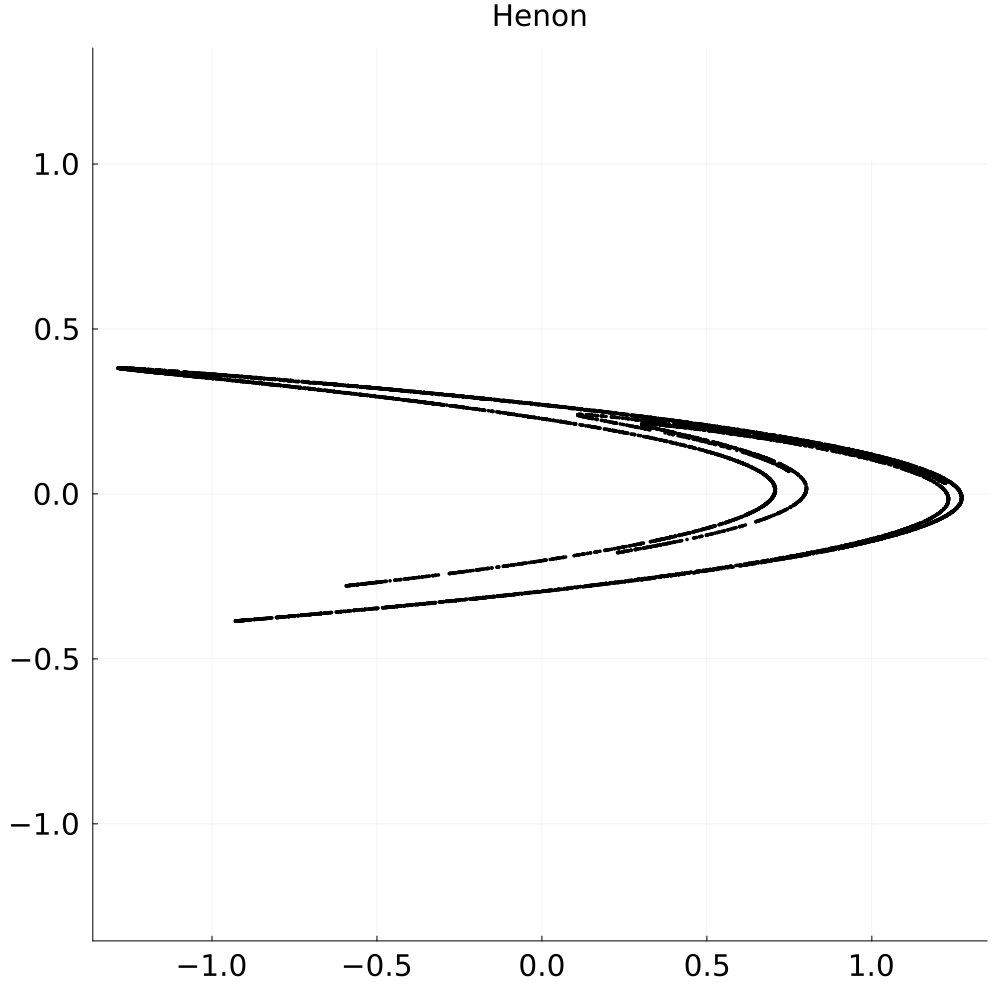}
	\includegraphics[width=0.35\linewidth, height=0.5\textheight, keepaspectratio]{\figs 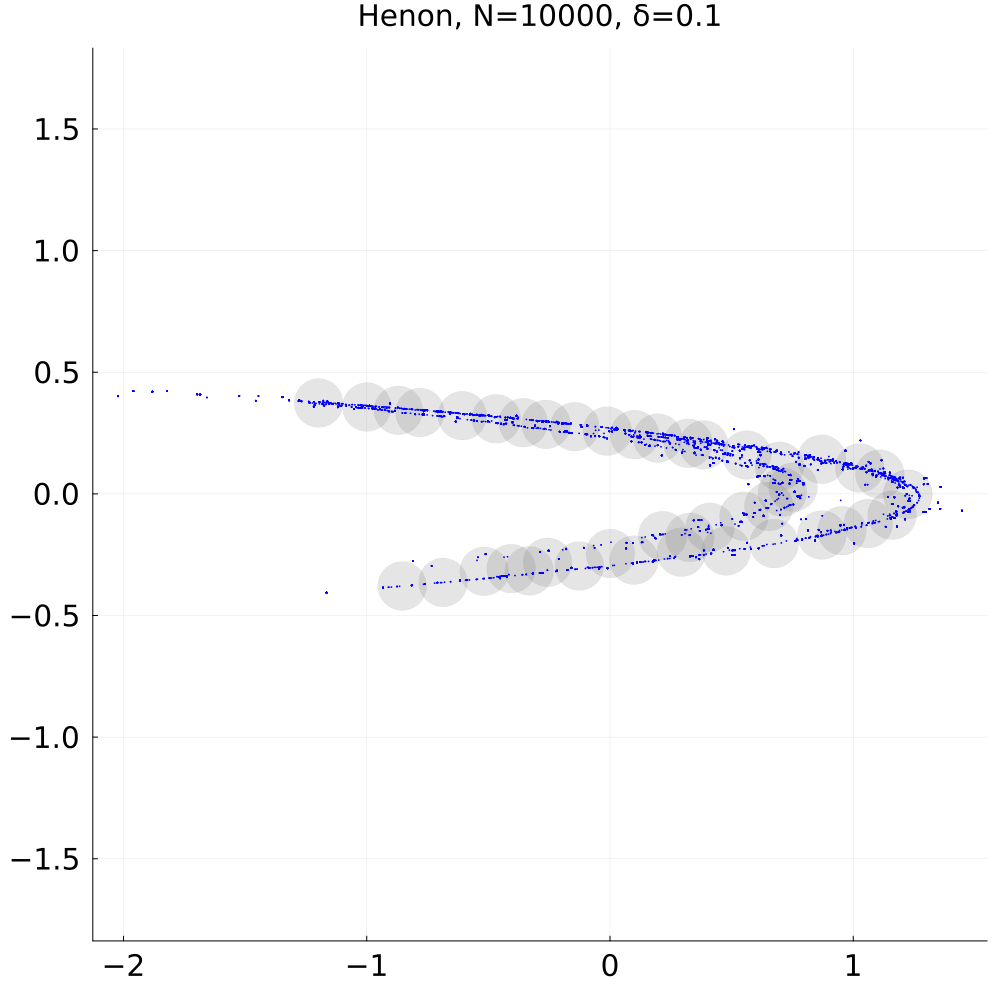}
	\includegraphics[width=0.28\linewidth, height=0.5\textheight, keepaspectratio]{\figs 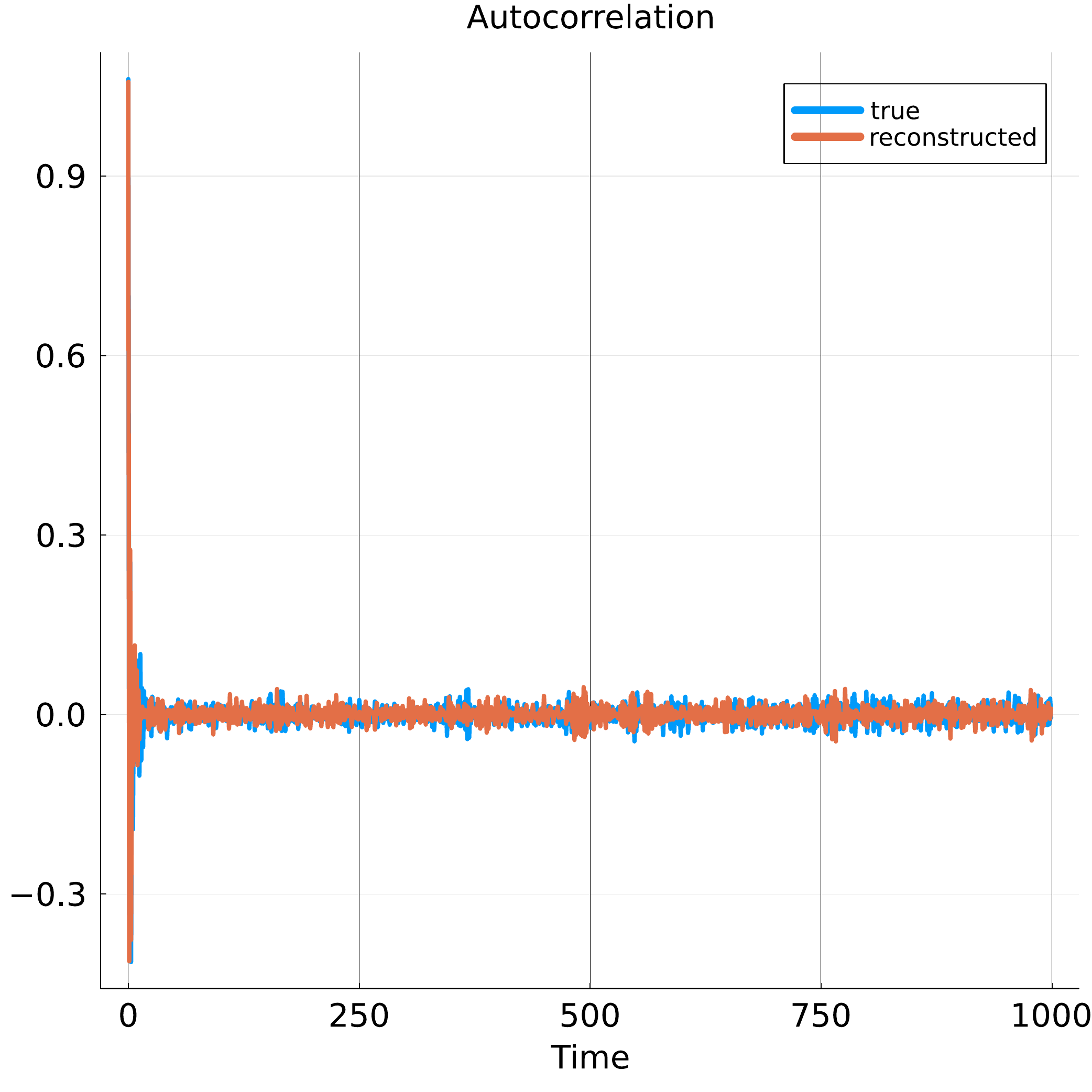}
	\caption{Reconstructing the Henon map \eqref{eqn:def:Henon}.}
	\label{fig:Henon}
\end{figure}

In this section we put to test our algorithms. The systems that we investigate are as follows : 
\begin{enumerate}
	\item \textbf{Lorenz 63.} This is the benchmark problem for studying chaotic phenomenon in three dimensions, which is the lowest dimension in which chaos is possible.
	\begin{equation} \label{eqn:vec:l63}
		\begin{split}
			\frac{d}{dt} x_1(t) &= \sigma(x_2-x_1) \\
			\frac{d}{dt} x_2(t) &= x_! (\rho-x_3) - x_2\\
			\frac{d}{dt} x_3(t) &= x_1 x_2 - \beta x_3
		\end{split}
	\end{equation}
	\item \textbf{Henon map.} The Henon map \cite{BenedicksCarleson1991Henon, AlvesEtAl2010Henon} is one of the simplest 2-dimensional dynamics that represent chaotic behavior.
	\begin{equation} \label{eqn:def:Henon}
		(x,y) \,\mapsto\, \paran{ 1+y-ax^2 \,,\, bx }
	\end{equation}
	\item \textbf{Lorenz 96.} The Lorenz 96 model is a dynamical system formulated by Edward Lorenz in 1996
	\begin{equation} \label{eqn:L96}
		\frac{d}{dt} x_n(t) = \paran{ x_{n+1} - x_{n-2} } x_{n-1} - x_n + F , \quad \forall n\in 1, \ldots, N .
	\end{equation}
	This model mimics the time evolution of an unspecified weather measurement collected at $N$ equidistant grid points along a latitude circle of the earth \cite{LorenzEmanuel1998opti}. The constant $F$ is known as the \emph{forcing} function. The system \eqref{eqn:L96} is benchmark problem in data assimilation.
\end{enumerate}

The numerical procedure outlined in Figure \ref{fig:outline1} was applied to all the systems. Some of the parameters involved in the procedure were kept constant, as declared in Table \ref{tab:param1}. All the experiments involved a variation in parameters $\delta$ and $N$. These experiments are listed in Table \ref{tab:param2}.

\paragraph{Bandwidth selection} We now present our method of selecting a bandwidth $\theta$ for constructing a kernel over a dataset $\calD$. It depends on the choice of a threshold $\eta \in (0,1)$, and a zero threshold $\theta_{zero} = 10^{-14}$. We first choose a subsample $\calD'$ of the data set $\calD$ and construct the set $\calS$ of all possible pairwise squared distances. Usually, it suffices to choose $\calD'$ to be 0.1 fraction of the dataset $\calD$, equidistributed throughout $\calD$. Now set another threshold : 
\begin{equation} \label{eqn:bw_select}
	\theta = -1/ \ln \paran{ \theta_{zero} }.
\end{equation}
Given the $\eta$ , we set $\theta$ such that $\eta$-fraction of the numbers in $\calS$ are less than $\theta$. In other words, with this choice of $\theta$, an $\eta$ fraction of the set $\SetDef{ e^{ - \norm{x-x'}^2 / \theta} }{x,x'\in \calD'}$ are greater than $\theta_{zero}$.

We next discuss various ways to evaluate the results.

\paragraph{Comparing topology} Our primary motivation was to recreate the dynamics as well as the invariant set $X$ supporting the measure $\mu$. Some care is needed to compare an invariant region of the reconstructed system $F$ with $h(X)$. They may not always be comparable as manifolds. One of the characteristic features of chaotic attractors \cite[e.g.]{Das_Yorke_PartialCtrl_2016, DasSaddles2015, DasYorke2020} is that they have a fractal geometry, which is self-similar at all scales, and almost nowhere differentiable. A general mathematical description of chaotic sets is elusive, and there is no reliable numerical method to characterize chaotic sets from data. We avoid this issue by comparing the two subsets by the Hausdorff metric $\dist_{\text{Hauss}}$ : 
\begin{equation} \label{eqn:def:Hauss_metric}
	\begin{split}
		\dist_{\to} \paran{X, Y} := \sup_{x\in X} \inf_{y\in Y} \dist(X, Y) ,\\
		\dist_{\text{Hauss}} \paran{ A, B } &:= \paran{ \dist_{\to} \paran{A, B} , \quad \dist_{\to} \paran{B, A} }
	\end{split}
\end{equation}
The Hausdorff metric is a proper metric on the collection of compact sets in any Hausdorff space. One can similarly define an $L^1$-metric 
\begin{equation} \label{eqn:def:Hauss_metric_L!}
	\begin{split}
		\dist_{L1, \to} \paran{X, Y} := \int_{x\in X} \inf_{y\in Y} \dist(X, Y) d\mu(x),\\
		\dist_{\text{Hauss}, L1} \paran{ A, B } &:= \paran{ \dist_{L1, \to} \paran{A, B} , \quad \dist_{L1, \to} \paran{B, A} }
	\end{split}
\end{equation}
The $L^1$ Hausdorff metric will be our choice of metric to evaluate the accuracy of our Markov approximation.

\paragraph{Comparing invariant measures} An invariant measure of a dynamical system bears many characteristics of the dynamics otherwise not available from just the topology of the invariant set. Invariant measures are the basis of the ergodic theorem, which essentially says that time average of a function along a typical trajectory must equal its space average with respect the invariant measure. In fact this phenomenon of convergence can provide a complete phase portrait of the multiple invariant sets co-existing with in $\calM$ \cite[e.g.]{DasJim2017_SuperC, DSSY2017_QQ}. Ergodic convergence is special case of a property called decay of correlations. Given two square integrable functions $\alpha, \beta : \calM \to \real^d$, their correlation is a function of time given by : 
\[ \Corel(n; \alpha, \beta) := \left\langle \alpha , U^n \beta \right\rangle_{L^2(\mu)} = \int_{x\in \calM} \left\langle \alpha(x) \beta \paran{ f^n x } \right\rangle_{\real^d} d\mu(x). \]
Suppose $\alpha, \beta$ are the same functions $\Phi$. Then the correlation function is named \emph{auto-correlation} :
\[ \AutoCorel(n; \Phi) := \Corel(n; \Phi, \Phi) \]
The decay of the correlation and auto-correlation functions with time $n$ indicates the level of mixing in the system. Mixing is the measure theoretic analog of chaos. Moreover, the collection of all correlation functions determine a construct called the \emph{spectral-measure} of the dynamics \cite[see]{DGJ_compactV_2018}. The spectral measure is an equivalent description of the dynamics using the language of operators and operator families. Thus in conclusion, autocorrelation functions form a good metric to compare the invariant measures of two different processes.

In a numerical setup we cannot compute the autocorrelation function for ever function $\Phi$. Instead we shall concentrate on the cases when $\Phi$ is the collection of coordinate projects. This autocorrelation function is easy to evaluate on signals. Given a signal $\vec{X} = \paran{x_1, \ldots, x_N}$, its autocorrelation can be approximated as
\begin{equation} \label{eqn:dp0j4}
	\Delta(\vec{x})(t) := \frac{1}{N-T} \sum_{n=1}^{N-T} \left\langle x_n - \bar{x}, x_{n+t} - \bar{x} \right\rangle_{\real^d}, \quad 1\leq t \leq T.
\end{equation}
The quantity in \eqref{eqn:dp0j4} will can be computed both for the original timeseries as well as the timeseries obtained by simulating the Markov process, for various lead-times $T$.

Table \ref{tab:param2} summarizes the results of our experiments. This ends the presentation of the theory, algorithm and experimental results. We next summarize the various aspects of the technique.

\section{Conclusions} \label{sec:conclus}

We have thus presented how an ergodic dynamical system $(\Omega, \mu, f)$ may be realized as the zero-noise limit of a Markov process. If the system is stochastically stable, then the Markov process also provides an approximation of its ergodic / statistical properties. These Markov processes can be implemented though data-driven algorithms. The experiments conducted with these algorithms confirmed the convergence claimed and proved in Theorem \ref{thm:1}. 

There are two main requirements in any data driven approach to dynamical system. Firstly, the methods must not rely on a fore-knowledge of the original phase space $\Omega$ or of the map $\mu$. The reconstruction must be entirely in a space created out of data, such as Euclidean spaces. Secondly the methods must preserve the invariant region and measures as seen from from observation maps. The loss of stability of the targeted invariant region is a common challenge faced by most numerical methods.

The procedure we have described fulfills both these requirements. Moreover, it has the following useful features :
\begin{enumerate}
\item The techniques are based on kernel methods. Kernel methods compute pairwise distances between data points, and convert the problem into an analysis of the kernel matrix. As a result, dimensionality issues are averted. One may have a low dimensional invariant set $X$ embedded in a high dimensional Euclidean space $\real^d$. The performance of kernel methods are not affected by the magnitude of $d$.
\item The technique places no strict requirement on the creation of the open covers. This allows many other possibilities, such as coverings using cubes and simplexes.
\item The main result Theorem \ref{thm:1} is stated in the context of the continuum measure $\mu$. This measure needs to be well sampled by data. For any choice of mesh size $\delta$, each cell $U_j$ of the cover must contain at least two data points. This simple requirement significantly reduces the demand for large data for larger and larger estimates.
\item One of the key ideas used in the technique was the use of Markov kernels to ensure that the ranges of the $\phi_{j\to i}$ remain confined to a convex hull. See their definition in \eqref{eqn:def:ph_ji}. This simple idea ensures the Hausdorff convergence of the support of the invariant measure for the stochastic process to the targeted invariant set $h(X)$.
\end{enumerate}	
	
	\bibliographystyle{unsrt_inline_url}
	\bibliography{\Path References,ref}
\end{document}